\documentclass{tac}

\usepackage{amssymb}
\usepackage{amsmath}
\usepackage{array}
\usepackage{bbm}
\usepackage{enumerate}
\usepackage[latin1]{inputenc}
\usepackage[T1]{fontenc}
\usepackage{shadow}
\usepackage[all,cmtip,2cell]{xy}
\UseAllTwocells
\usepackage{graphicx}
\usepackage{graphbox}
\usepackage{caption}
\usepackage{subcaption}
\usepackage{float}
\usepackage{multicol}
\usepackage{tikz}
\usepackage{tikz-cd}
\usepackage{longtable}
\usepackage{color}
\usetikzlibrary{matrix,arrows,decorations.pathmorphing}
\usepackage{pdfpages}

\usepackage[colorlinks=true]{hyperref}
\hypersetup{allcolors=[rgb]{0.1,0.1,0.4}}

\newcommand{\CC}{\mathcal{C}}
\newcommand{\CCC}{\mathfrak{C}}
\newcommand{\DDD}{\mathfrak{D}}
\newcommand{\EE}{\mathcal{E}}

\newcommand{\OO}{\mathcal{O}}

\newcommand{\DD}{\mathcal{D}}

\newcommand{\G}{\mathbb{G}}
\newcommand{\GG}{\mathcal{G}}
\newcommand{\N}{\mathbb{N}}
\newcommand{\sd}{D^s}

\newcommand{\T}{\mathcal{T}}

\newcommand{\quot}{\delimiter"502F30E\mathopen{}}

\mathrmdef{Fun}
\mathrmdef{Hom}
\mathrmdef{gSet}
\mathrmdef{Set}
\mathrmdef{id}
\mathrmdef{Id}
\mathrmdef{Map}
\mathrmdef{Adj}
\mathrmdef{Alg}
\mathrmdef{Comp}
\mathrmdef{Cat}
\mathrmdef{pd}
\mathrmdef{C}
\mathrmdef{Sesq}
\mathrmdef{Tree}
\mathrmdef{im}
\mathrmdef{obj}
\mathrmdef{sk}

\definecolor{myblack}{RGB}{0,0,0}
\definecolor{myred}{RGB}{255,0,0}
\definecolor{myblue}{RGB}{0,0,255}
\definecolor{mygreen}{RGB}{0,120,0}
\definecolor{mypurple}{RGB}{120,0,120}
\definecolor{myorange}{RGB}{255,120,0}
\definecolor{mylightblue}{RGB}{0,255,255}
\definecolor{myolive}{RGB}{120,120,0}
\definecolor{myteal}{RGB}{0,120,120}

\title{Simple string diagrams and $n$-sesquicategories}  

\author{Manuel Ara\'{u}jo}     
\address{Department of Computer Science and Technology, University of Cambridge}
\eaddress{manuel.araujo@cst.cam.ac.uk}


\copyrightyear{2021}

\thanks{I would like to thank Richard Garner for suggesting the inductive characterization of $n$-sesquicategories and the anonymous referee for many useful comments and suggestions. This work was partially supported by the FCT project grant \textbf{Higher Structures and Applications}, PTDC/MAT-PUR/31089/2017. }

\begin{document}

\maketitle                 

\begin{abstract}	

We define a monad $T_n^{\sd}$ whose operations are encoded by simple string diagrams and we define $n$-sesquicategories as algebras over this monad. This monad encodes the compositional structure of $n$-dimensional string diagrams. We give a generators and relations description of $T_n^{\sd}$, which allows us to describe $n$-sesquicategories as globular sets equipped with associative and unital composition and whiskering operations. One can also see them as strict $n$-categories without interchange laws. Finally we give an inductive characterization of $n$-sesquicategories.

\end{abstract}

\section{Introduction}

We are interested in developing a theory of $n$-categories based on string diagrams. In this paper we define a monad $T_n^{\sd}$ which encodes the compositional structure of $n$-dimensional string diagrams. We call $T_n^{\sd}$-algebras \textbf{$n$-sesquicategories}. For $n=2$ these are already known as sesquicategories (see \cite{street_sesquicat}). We also give a finite presentation of this monad by generators and relations, and we give an inductive description of $n$-sesquicategories as categories equipped with a lift of the $\Hom$ functor to the category of $(n-1)$-sesquicategories.

\subsection{Results}

We define an $n$-globular set $\sd_n$ such that $\sd_n(k)$ is the set of $k$-dimensional simple string diagrams. These diagrams play a role analogous to that of globular pasting diagrams in the definition of strict $n$-categories (see \cite{OperadsCats}). We then define the monad $T_n^{\sd}$ on $n$-globular sets by the formula $$T_n^{\sd}(X)(m):=\coprod_{D\in\sd_n(m)}\gSet_n(\widehat{D},X).$$ This means that an $n$-sesquicategory is an $n$-globular set $\CC$ together with a choice of composite for each $\CC$-labeled simple string diagram. We give a finite presentation of this monad by generators and relations, which allows us to give the following description of $n$-sesquicategories. Below, we denote by $\CC_i$ the set of $i$-morphisms in $\CC$ and we write $|x|=i$ when $x$ is an $i$-morphism.

\begin{theorem}\label{genrel}
	
An $n$-sesquicategory is an $n$-globular set $\CC$ equipped with 

\begin{itemize}
	
	\item binary operations $\circ^{\CC}_{i,j}\colon \CC_{i}\times_{\CC_{m-1}}\CC_j\to\CC_M$ for $1\leq i,j\leq n$, where $m=\min\{i,j\}$, $M=\max\{i,j\}$ and the pullback is taken with respect to the maps $s^{i-m+1}$ and $t^{j-m+1}$;
	\item unary operations $u_i^{\CC}\colon \CC_{i-1}\to\CC_i$ for $1\leq i\leq n$;
	
\end{itemize} subject to equations

\begin{itemize}
	
	\item $s(x\circ y)=s(x)\circ y$ and $t(x\circ y)=t(x)\circ y$, if $|x|>|y|$;
	
	\item $s(x\circ y)=x\circ s(y)$ and $t(x\circ y)=x\circ t(y)$, if $|x|<|y|$;
	
	\item $s(x\circ y)=s(y)$ and $t(x\circ y)=t(x)$, if $|x|=|y|$;
	
	\item $s(u(x))=t(u(x))=x$;
	
	\item $(x\circ y)\circ z=x\circ (y\circ z)$ when the two smallest elements in $\{|x|,|y|,|z|\}$ are equal;
	
	\item $x\circ (y\circ z)=(x\circ y)\circ(x\circ z)$ when $|x|<|y|,|z|$;
	
	\item $(x\circ y)\circ z=(x\circ z)\circ(y\circ z)$	when $|z|<|x|,|y|$;
	
	\item $u(x)\circ y = y$ when $|x|<|y|$;
	
    \item $u(x)\circ y = u(x\circ y)$ when $|x|\geq |y|$;
	
	\item $x\circ u(y) = x$ when $|x|>|y|$;
	
	\item $x\circ u(y)=u(x\circ y)$ when $|x|\leq |y|$;
	
\end{itemize}
	
\end{theorem}	 

Moreover, we give an inductive description of $n$-sesquicategories, analogous to the inductive definition of strict $n$-categories as categories enriched in $\Cat_{n-1}$. Denote the category of $T_n^{\sd}$-algebras and algebra homomorphisms by $\Sesq_n$.

\begin{theorem}
	
A $1$-sesquicategory is simply a category. An $n$-sesquicategory is a category together with a lift of its $\Hom$ functor	\[\xymatrix{ & \Sesq_{n-1}\ar@{}[d]|-{=}\ar[rd]^{(-)_0} &  \\\CC^{op}\times\CC\ar@{.>}[ru]^{\underline{\Hom}_\CC}\ar[rr]_{\Hom_\CC} & & \Set .}\]
	
\end{theorem}	

\subsection{Future work}

\textbf{String diagrams} are a useful computational tool in many kinds of higher algebraic settings. In \cite{thesis}, \cite{paper1}, \cite{paper2} and \cite{adj4} we have used string diagrams in strict $3$ and $4$-categories to prove results about fibrations of mapping groupoids and coherence for adjunctions. We would like to be able to extend this method beyond strict $n$-categories. To this effect, we want to define a notion of \textbf{semistrict $n$-category} which is more general than a strict $n$-category, but which still admits composites for string diagrams. Our strategy is to have string diagrams built into the definition, so that a semistrict $n$-category would be defined as an $n$-globular set $\CC$ equipped with a way of composing string diagrams labeled in $\CC$.

Since string diagrams are composed simply by juxtaposition, it is natural to require that our $n$-categories should have strictly associative and unital composition operations. However, unlike for strict $n$-categories, we don't require interchange laws to hold: $$\includegraphics[scale=1.5,align=c]{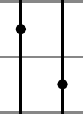}\neq \includegraphics[scale=1.5,align=c]{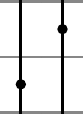}.$$ We also don't allow certain forms of composition which would cause difficulties when drawing higher dimensional diagrams. As an example, the following diagram is not allowed:$$\includegraphics[scale=1.5,align=c]{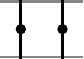}$$ These considerations lead to our definition of the monad $T_n^{\sd}$.

The algebras over $T_n^{\sd}$ are not yet our desired notion of semistrict $n$-category. To obtain semistrict $n$-categories, one would need to add coherence cells implementing weak versions of the equations which are satisfied in strict $n$-categories, obtaining a monad $T_n^{ss}$. This means that whenever two string diagrams correspond to the same pasting diagram, they should be connected by a composite of coherence cells. In an upcoming paper, we do this for $n=3$. Namely, we define a monad $T_3^{ss}$ by adding to $T_3^{\sd}$ an operation encoding the interchanger of $2$-morphisms depicted below, satisfying cancellation and Yang-Baxter equations. $$\includegraphics[scale=1.5,align=c]{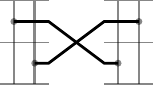}$$ We call $T_3^{ss}$-algebras semistrict $3$-categories and show that they are equivalent to Gray categories. This gives a new definition of Gray categories, based entirely on string diagrams. We are also working on extending this procedure to higher dimensions. 

\subsection{Relation to \cite{Jamie}}

 In \cite{Jamie}, the authors introduce the concept of an \textbf{$n$-signature}, which is a set of generating morphisms in each dimension $k\leq n$. These generating morphisms can be composed in certain ways to form \textbf{$n$-diagrams}, which are $n$-dimensional string diagrams. They define an \textbf{$n$-sesquicategory} as an $(n+1)$-signature. The authors then introduce the notion of \textbf{homotopy generators}, which implement certain forms of coherence which are expected in $n$-categories. They propose lists of homotopy generators in dimensions $\leq 5$ and then conjecture that there should be some notion of a \textbf{semistrict $n$-category} for which semistrict $n$-categories are presented by $(n+1)$-signatures with homotopy generators. They claim that for $n=1,2,3$ these are categories, strict $2$-categories and Gray categories, respectively. They also claim that for $n=4$ these should correspond to Crans' $4$-tas.
 
 In this paper, we start the work of finding this notion of semistrict $n$-category. An $(n+1)$-signature corresponds in our language to an $(n+1)$-computad for $T_n^{\sd}$ (a computad for $T$ is just a presentation of a $T$-algebra by generators and relations, see \cite{CSPThesis}). So what is called an $n$-sesquicategory in \cite{Jamie} is for us a $(n+1)$-computad for $T_n^{\sd}$, which means a presentation for an $n$-sesquicategory.  Once we define the monad $T_3^{ss}$, a $4$-signature with homotopy generators will correspond to a $4$-computad for $T_3^{ss}$. We are working on extending this to higher dimensions, where an $(n+1)$-signature with homotopy generators would correspond to an $(n+1)$-computad for $T_n^{ss}$. The advantage of having the monads $T_n^{ss}$ would be that we can talk about semistrict $n$-categories without choosing finite presentations. This is what allows us to prove in an upcoming paper that semistrict $3$-categories are the same thing as Gray categories, as predicted in \cite{Jamie}.
 
 \subsection{Relation to other work}
 
The use of string diagrams in the context of tensor categories and $2$-categories of various kinds is by now well established (see \cite{nlabstring} for a survey). There has also been work on extending string diagrams to dimension three, viewing string diagrams as certain kinds of embedded manifolds (\cite{bms}, \cite{hummon}, \cite{trimble}). Here we define string diagrams as combinatorial structures, which has the advantage of allowing for finite, computable descriptions and an algebraic definition of semistrict $n$-categories. 

The problem of relating topological and algebraic concepts in higher category theory is  also recently discussed in \cite{fct} and \cite{mfd_diag}. One of the authors has also proposed a definition of associative $n$-category in his PhD Thesis \cite{dorn}, which seems to be closely related to what we have called semistrict n-category here. An online proof assistant for diagrammatic calculus in associative $n$-categories has also appeared (see \cite{rv} and \cite{homotopyio}).

Finally, we note that in the preprint \cite[Section 8]{makkai} there is a description of strict $\omega$-catgories as globular sets equipped with binary and unary operations, satisfying certain laws which are exactly analogous to the ones in Theorem \ref{genrel}, plus an additional law expressing the Godement interchange rules.

During the referee process, we became aware of the paper \cite{rewriting_gray}, where a notion of \textbf{$n$-precategory} is introduced, which is equivalent (by Theorem \ref{genrel}) to our notion of $n$-sesquicategory.

In \cite{sesqui_comp} we prove that the category of computads for the monad $T_n^{\sd}$ is a presehaf category and we describe a string diagram notation for morphisms in an $n$-sesquicategory generated by a computad.

\thirdleveltheorems

\section{Simple string diagrams}

In this section we define simple string diagrams and the monad $T_n^{D^s}$ whose operations are encoded by these diagrams and whose algebras we name $n$-sesquicategories.

We start by recalling the definition of strict $n$-categories as algebras over a monad $T_n$ whose operations are encoded by globular pasting diagrams. Then we define simple string diagrams and explore their relation to globular pasting diagrams. Finally, we define the monad $T_n^{D^s}$ by the same formula as $T_n$, with globular pasting diagrams replaced by simple string diagrams. 

\subsection{Globular pasting diagrams and strict $n$-categories}

Let $\G_n$ be the $n$-globe category, with set of objects $\{0,\cdots,n\}$ and morphism sets generated by $s,t\colon k\to k+1$ satisfying the usual globularity relations $ts=ss$ and $tt=st$. Let $$\gSet_n:=\Fun(\G_n^{op},\Set)$$ the category of \textbf{$n$-globular sets}. So an $n$-globular set consists of sets $X(k)$ of \textbf{$k$-cells} for $0\leq k\leq n$, together with \textbf{source} and \textbf{target} maps $s,t\colon X(k)\to X(k-1)$ such that $st=ss$ and $tt=ts$. We sometimes denote $X(k)$ by $X_k$ when this does not cause confusion with other subscripts. We also write $|x|=k$ when $x\in X(k)$. We have adjoint functors $\xymatrix{\gSet_k\ar@<1ex>[r]_{\perp} & \gSet_n\ar@<1ex>[l]}$ for $k\leq n$, given by restriction and extension by $\emptyset$. We sometimes use these implicitly, writing $\gSet_n(X,Y)$ or $\gSet_k(X,Y)$ when $X\in\gSet_k$ and $Y\in\gSet_n$, for example.

We explain the definition of a strict $n$-category as an algebra over a certain monad on $n$-globular sets, following \cite[Chapter 8]{OperadsCats}.

\begin{notation}

Denote by $\Delta_a$ the augmented simplex category, consisting of finite totally ordered sets and order preseving maps. Denote $\langle n \rangle :=\{1\leq\cdots\leq n\}$ for $n\in\N$.
 
\end{notation}

\begin{definition} 

A \textbf{$k$-stage level tree} is a diagram $$\xymatrix@1{X_k\ar[r] & \cdots \ar[r] & X_1\ar[r] & \langle 1 \rangle}$$ in $\Delta_a$. An \textbf{isomorphism} of $k$-stage level trees is an isomorphism of diagrams. An \textbf{unordered isomorphism} of $k$-stage level trees is an isomorphism of the underlying diagrams of sets. Denote by $\T_k$ and $\T^u_k$ the groupoids of $k$-stage level trees with isomorphisms and unordered isomorphisms, respectively.

\end{definition}

\begin{definition}

A \textbf{globular $k$-pasting diagram} $\pi$ is a $k$-stage level tree of the form \[\xymatrixcolsep{2.5pc}\xymatrix@1{\langle\ell_k(\pi)\rangle\ar[r]^-{v^k(\pi)} & \cdots \ar[r]^-{v^2(\pi)} & \langle\ell_1(\pi)\rangle\ar[r]^-{v^1(\pi)} & \langle\ell_0(\pi)\rangle=\langle 1 \rangle}.\] Its source $s(\pi)$ and target $t(\pi)$ are equal and given by truncation. This defines an $n$-globular set which we call $\pd_n$. We sometimes denote $v^m(\pi)$ and $\ell_m(\pi)$ simply by $v^m$ and $\ell_m$.
	
\end{definition}

\begin{remark}

Any $k$-stage level tree $X$ is uniquely isomorphic to a a unique globular $k$-pasting diagram.

\end{remark}

The set $\pd(1)$ of globular $1$-pasting diagrams is just the set of natural numbers. The pictorial representation of the $1$-pasting diagram corresponding to $m$ is a string of $m$ composable arrows $$\xymatrix{\bullet\ar[r] & \bullet\ar[r] & \cdots\ar[r] & \bullet}.$$ 

The $2$-stage level tree on the left corresponds to the pasting diagram on the right. \begin{center}\begin{tabular}{lcccr}\includegraphics[scale=2,align=c]{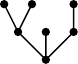} & & & & $\xymatrix@1{\bullet\ruppertwocell\rlowertwocell\ar[r] & \bullet \ar[r] &\bullet \rtwocell & \bullet}$\end{tabular}\end{center}

To each $k$-pasting diagram $\pi\in\pd_n(k)$, it is possible to associate a $k$-globular set $\widehat{\pi}\in\gSet_k$ whose cells are the the ones sugested by its pictorial representation.  If we take $\pi$ to be the $2$-dimensional example above, then $\widehat{\pi}$ has four $0$-cells, six $1$-cells and three $2$-cells, and the source and target maps can be easily read from the picture. A precise inductive definition can be found in \cite{OperadsCats}. Here we give a different definition, which is easily shown to be equivalent.

\begin{notation}
	
For $k\in\N$, denote by $\C_k$ the $k$-globular set $\Hom_{\G_k}(-,k)$ represented by $k$.
	
\end{notation}

\begin{definition}
	
Let $\pi$ be a globular $k$-pasting diagram. We define the $k$-globular set $\widehat{\pi}$ by a sequence of cell attachments, meaning we inductively define $m$-globular sets $(\widehat{\pi})^m_i$, with $0\leq m \leq k$ and $0\leq i \leq \ell_i$, and take $$\widehat{\pi}:=(\widehat{\pi})^k_{\ell_k}.$$ We define $(\widehat{\pi})^0_0=\emptyset$, $(\widehat{\pi})^0_1=C_0$ and $(\widehat{\pi})^m_0=(\widehat{\pi})^{m-1}_{\ell_{m-1}}$. For $m\geq 1$ we build $(\widehat{\pi})^m_i$ from $(\widehat{\pi})^m_{i-1}$ by attaching an $m$-cell $x^m_i$ as follows. 

If $i=\min\{j\in\pi(m):v^m(j)=v^m(i)\}$ then the attachment is described by the following pushout of $m$-globular sets:

\[\xymatrix@1{C_{m-1}\ar[r]^s \ar[d]_{x^{m-1}_{v^m(i)}} & C_m \ar[d]^{x^m_i} \\ (\widehat{\pi})^m_{i-1}\ar[r] & (\widehat{\pi})^m_{i}} \] Otherwise it is described by the following pushout of $m$-globular sets:

\[\xymatrix@1{C_{m-1}\ar[r]^s \ar[d]_{t(x^m_{i-1})} & C_m \ar[d]^{x^m_i} \\ (\widehat{\pi})^m_{i-1}\ar[r] & (\widehat{\pi})^m_{i}} \]

\end{definition}

This definition corresponds to the simple procedure of moving up the tree, attaching an $m$-cell for each node at level $m$ by gluing its source to a prescribed $(m-1)$-cell. If the node in question is the first one to be connected to its parent $(m-1)$-node, then its source is attached to the $(m-1)$-cell corresponding to this node. Otherwise, it is attached to the target of the $m$-cell corresponding to the previous node at level $m$ with the same parent node.

\begin{definition}
	
	Given a simple $k$-pasting diagram $\pi$ and a $k$-globular set $X$, we define an \textbf{$X$-labelling} of $\pi$ to be a map of $k$-globular sets $\hat{\pi}\to X$.
	
\end{definition}

This corresponds to the idea of putting labels on each dot and arrow in the graphical depiction of the pasting diagram, with each label being a cell in $X$ of the correct dimension, satisfying source and target compatibility. 

Consider the functor

$$T_{n}\colon \gSet_n\to\gSet_n$$ where $T_{n}(X)$ has $$T_{n}(X)_{k}=\coprod_{\pi\in\pd_n(k)}\gSet_n(\hat{\pi},X)$$ as its set of $k$-cells. This functor can be given a monad structure, by interpreting a pasting diagram labelled by pasting diagrams as a pasting diagram.

\begin{definition}
	
	An \textbf{$n$-category} is an algebra over the monad $T_n$. A \textbf{functor} between $n$-categories is a morphism of algebras over $T_n$. We denote the category $\Alg_{T_n}$ of $T_n$-algebras and $T_n$-algebra morphisms by $\Cat_n$. 
	
\end{definition} 

This agrees with the standard definitions of strict $n$-categories and strict functors between them.

\subsection{Simple string diagrams}

Now we define simple string diagrams and introduce their pictorial representations.

\begin{definition} A \textbf{$k$-stage level tree with crossings} is a diagram $$\xymatrix@1{X_k\ar[r] & \cdots \ar[r] & X_1\ar[r] & \langle 1 \rangle}$$ of finite totally ordered sets and maps which are not required to be order preserving. An \textbf{isomorphism} of $k$-stage level trees with crossings is a levelwise order preserving isomorphism of diagrams. An \textbf{unordered isomorphism} of $k$-stage level trees with crossings is an isomorphism of the underlying diagrams of sets. We denote the corresponding groupoids by $\T_k^{\chi}$ and $\T_k^{\chi, u}$.

\end{definition}

\begin{definition}
	A \textbf{simple $k$-string diagram} $D$ is a $k$-stage level tree with crossings of the form \[\xymatrixcolsep{2.5pc}\xymatrix@1{\langle\ell_k(D)\rangle\ar[r]^-{v^k(D)} & \cdots \ar[r]^-{v^2(D)} & \langle\ell_1(D)\rangle\ar[r]^-{v^1(D)} & \langle\ell_0(D)\rangle=\langle 1 \rangle}.\] The source $s(D)$ and target $t(D)$ of $D$ are equal and given by truncation. This defines an $n$-globular set, which we denote by $\sd_n$. We sometimes denote $v^m(D)$ and $\ell_m(D)$ simply by $v^m$ and $\ell_m$.
	
\end{definition}

\begin{remark}

Any $k$-stage level tree with crossings $X$ is uniquely isomorphic to a a unique simple $k$-string diagram.

\end{remark}

\begin{notation} Each $v^i(D)$ is simply a sequence of length $\ell_i(D)$ taking values in $\langle\ell_{i-1}(D)\rangle$ and we usually denote a diagram by the corresponding sequence of sequences. We usually replace the sequence $v^1(D)$ by the natural number $\ell_1(D)$. 
\end{notation}

\begin{example} Consider the simple $2$-string diagram $\langle3\rangle\to\langle3\rangle\to\langle1\rangle$ where the map $\langle3\rangle\to\langle3\rangle$ sends $1,3\mapsto 1$ and $2\mapsto 3$. This $2$-stage level tree with crossings can be depicted as \begin{center}\includegraphics[scale=2,align=c]{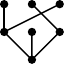}.\end{center} We can also denote this simple $2$-string diagram by the sequence of sequences $((1,1,1),(1,3,1))$ or simply by $(3,(1,3,1))$.

\end{example} 

We now explain how to think about these trees as string diagrams by giving another graphical representation. We always read $k$-diagrams from left to right when $k$ is odd and from top to bottom when $k$ is even. The map $\ell_1\colon \sd(1)\to\N$ is a bijection and we interpret the element in $\sd(1)$ corresponding to the natural number $m\in\N$ as a string of $m$ composable morphisms, which we depict as $m$ dots on a line. For example, for $m=2$ we have $$\includegraphics[scale=1.5]{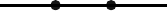}.$$

Now consider $D=(m,(k))\in\sd(2)$, for some $k\leq m$. We interpret this as a string of $m$ horizontally composable $2$-morphisms, with all being the identity $2$-morphism except the $k^{th}$ one. For example, we depict the simple $2$-string diagram $(5,(2))$ as $$\includegraphics[scale=1.5]{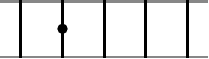}.$$ A generic element $D\in\sd(2)$ is then just a sequence of these blocks, which stack vertically from top to bottom. For example, here is a picture of $(3,(1,3,1)):$ $$\includegraphics[scale=1.5]{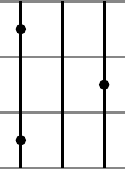}.$$

It should now be clear how the process is repeated in higher dimensions. Here is a picture of the simple $3$-string diagram $(3,(2,1),(1,1))$: \begin{center}\includegraphics[scale=1.5]{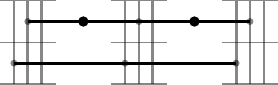}.\end{center} And here is a picture of the simple $4$-string diagram $(2,(1,2),(2,1),(2))$: $$\includegraphics[scale=1.5]{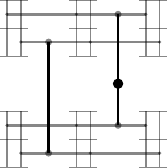}.$$

\subsection{Relation to globular pasting diagrams}

We now define a map of $n$-globular sets $$\pi\colon \sd_n\to\pd_n$$ and establish some of its properties.

\begin{notation}
	
When $X$ and $Y$ are totally ordered sets, we denote by $X\times^l Y$ their product equipped with the lexicographical order. This order is defined by declaring $(x_1,y_1)\leq (x_2,y_2)$ whenever one of the following conditions is satisfied:

\begin{itemize}

 \item $x_1<x_2$; 
 \item $x_1=x_2$ and $y_1\leq y_2$.
 
 \end{itemize}
	
\end{notation}	

\begin{definition}
	
We define a functor $\overline{(-)}:\T_k^{\chi,u}\to\T_k^u$ as follows. Given a  $k$-stage level tree with crossings \[X = \left( \xymatrix@1{X_k\ar[r]^-{v^k} & \cdots \ar[r]^{v^2} & X_1\ar[r]^-{v^1} & \langle 1 \rangle}\right) ,\] we define a $k$-stage level tree \[\overline{X}:= \left( \xymatrix@1{\overline{X}_k\ar[r]^-{v^k} & \cdots \ar[r]^{v^2} & \overline{X}_1\ar[r]^-{v^1} & \langle 1 \rangle}\right) ,\] where $\overline{X}_i$ denotes a total ordering on the underlying set of $X_i$ which makes the maps $v^i$ order preserving. We define this ordering inductively, by setting $\overline{X}_1=X_1$ and taking $\overline{X}_{i+1}$ to be the unique total order such that the map \[\xymatrixcolsep{4pc}\xymatrix@1{\overline{X}_{i+1}\ar[r]^-{v^{i+1}\times\id} &  \overline{X}_i \times^l X_{i+1}}\] is order preserving. Given an unordered isomorphism $\sigma\colon X\to Y$ of $k$-stage level trees with crossings, we define $\overline{\sigma}\colon \overline{X}\to\overline{Y}$ to be the unordered isomorphism of $k$-stage level trees with $\overline{\sigma}_i=\sigma_i$ as maps of underlying sets. 
	
\end{definition}

\begin{definition}

Given a simple $k$-string diagram $D$ we write $\pi(D)$ for the unique globular $k$-pasting diagram isomorphic to $\overline{D}$. This commutes with taking sources and targets and so defines a map of $n$-globular sets $\pi\colon \sd_n\to\pd_n$.
 
\end{definition}

\begin{example}
	
We have $$\pi\left( \includegraphics[scale=1.5,align=c]{stringdiag/2d_2.pdf} \right)=\xymatrix@1{\bullet\ruppertwocell\rlowertwocell\ar[r] & \bullet \ar[r] &\bullet \rtwocell & \bullet}.$$ This is because $\overline{(-)}$ acts by removing crossings, so it maps \begin{tabular}{lcr}\includegraphics[scale=2,align=c]{pasting/tree2.pdf} & $\mapsto$ & \includegraphics[scale=2,align=c]{pasting/tree.pdf}  \end{tabular} .
	
\end{example}

\begin{lemma}
	
	The map $\pi\colon \sd_n\to\pd_n$ is surjective. 
	
\end{lemma}	

\begin{proof}

Any $k$-level tree is also a $k$-level tree with crossings, which just happens to not have any crossings.\end{proof}

Now we turn to the question of identifying all simple $k$-string diagrams which are mapped to a particular globular $k$-pasting diagram. 

\begin{lemma}
	
Let $\sigma\in\T_k^{\chi,u}(X,Y)$. Then $\overline{\sigma}\in\T_k(\overline{X},\overline{Y})\subset\T_k^u(\overline{X},\overline{Y})$ if and only if each $\sigma_i$ is order preserving on the fibres of $v^i$.
		
\end{lemma}	

\begin{proof}
	
We have $X_1=\overline{X}_1$ and $Y_1=\overline{Y}_1$, and the fibre of $v^1$ is all of $X_1$, so there is nothing to prove for $i=1$.
Now consider the following diagram: \[\xymatrixcolsep{4pc}\xymatrix@1{\overline{X}_{m+1} \ar[r]^-{v^{m+1}\times \id}\ar[d]_{\overline{\sigma}_{m+1}} & \overline{X}_{m}\times_{X_m}^l X_{m+1}\ar[d]^{\overline{\sigma}_m\times\sigma_{m+1}}  \\ \overline{Y}_{m+1} \ar[r]^-{u^{m+1}\times \id} & \overline{Y}_{m}\times_{Y_m}^l Y_{m+1}  }\] By definition of the lexicographical order, and because $\sigma_m$ is injective, $\overline{\sigma}_m\times\sigma_{m+1}$ is order preserving if and only if $\overline{\sigma}_m$ is order preserving and $\sigma_{m+1}$ is order preserving on the fibres of $v^m$. On the other hand, since the diagram commutes and the horizontal maps are order preserving isomorphisms, we know that $\overline{\sigma}_{m+1}$ is order preserving if and only if $\overline{\sigma}_m\times\sigma_{m+1}$ is order preserving. We conclude that $\overline{\sigma}_{m+1}$ is order preserving if and only if $\overline{\sigma}_m$ is order preserving and $\sigma_{m+1}$ is order preserving on the fibres of $v^m$. This allows us to prove the result by induction.\end{proof}	
		
The following is just a restatement of this Lemma.

\begin{proposition}
	
Let $C, D$ be simple $k$-string diagrams. Then $\pi(C)=\pi(D)$ if and only if $\ell_i(C)=\ell_i(D)$ and there exist permutations $\sigma_i\in\Sigma_{\ell_i(D)}$ which are order preserving on the fibres of $v^i(D)$, such that $v^i(C)=\sigma_{i-1}\circ v^i(D)\circ \sigma_i^{-1}$ for all $i$.	
	
\end{proposition}

\begin{example}
	
There are exactly three simple $2$-string diagrams, pictured on the left, which map to the globular $2$-pasting diagram on the right. \begin{center}\begin{tabular}{lcr}\includegraphics[scale=1.5,align=c]{stringdiag/2d_2.pdf}, \includegraphics[scale=1.5,align=c]{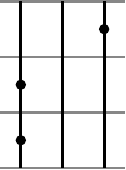}, \includegraphics[scale=1.5,align=c]{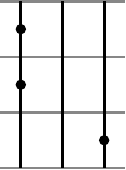} & $\mapsto$ & $\xymatrix@1{\bullet\ruppertwocell\rlowertwocell\ar[r] & \bullet \ar[r] &\bullet \rtwocell & \bullet}$ \end{tabular}\end{center}  These string diagrams are related to each other by permuting the heights of the three dots, but preserving the order of the two dots that appear on the same string.

\end{example}						

\subsection{Labelled simple string diagrams}

We now introduce labellings of simple string diagrams.

\begin{definition}
	
Let $D$ be a simple $k$-string diagram. We define the $k$-globular set $\widehat{D}$ by a sequence of cell attachments, meaning we inductively define $m$-globular sets $(\widehat{D})^m_i$, with $0\leq m \leq k$ and $0\leq i \leq \ell_m$, and take $$\widehat{D}:=(\widehat{D})^k_{\ell_k}.$$ We define $(\widehat{D})^0_0=\emptyset$, $(\widehat{D})^0_1=\C_0$ and $(\widehat{D})^m_0=(\widehat{D})^{m-1}_{\ell_{m-1}}$. For $m\geq 1$, we obtain $(\widehat{D})^m_i$ from $(\widehat{D})^m_{i-1}$ by attaching an $m$-cell $x^m_i$ as follows. 

If $i=\min\{j:v^m_j=v^m_i\}$ then the attachment is described by the following pushout of $m$-globular sets:

\[\xymatrix@1{C_{m-1}\ar[r]^s \ar[d]_{x^{m-1}_{v^m_i}} & C_m \ar[d]^{x^m_i} \\ (\widehat{D})^m_{i-1}\ar[r] & (\widehat{D})^m_{i}} \] Otherwise, we let $p=\max\{j<i:v^m_j=v^m_i\}$ and the attachment is described by the following pushout of $m$-globular sets:

\[\xymatrix@1{C_{m-1}\ar[r]^s \ar[d]_{t(x^m_{p})} & C_m \ar[d]^{x^m_i} \\ (\widehat{D})^m_{i-1}\ar[r] & (\widehat{D})^k_{i}} \]

\end{definition}	

It is clear from the above definition that $\widehat{\pi(D)}=\widehat{D}$.

\begin{definition}
	
Given a simple $k$-string diagram $D$, we define inclusions $\widehat{s(D)}\hookrightarrow\widehat{D}$ and $\widehat{t(D)}\hookrightarrow\widehat{D}$	as follows. Let $\widehat{s(D)}=\widehat{t(D)}$ be composed by attaching cells $$\{x^m_i:1\leq m\leq k-1\text{ and }1\leq i \leq \ell_m(D)\}.$$ and let $\widehat{D}$ be composed by attaching cells $$\{y^m_i:1\leq m\leq k\text{ and }1\leq i \leq \ell_m(D)\}.$$ Then the inclusion map $\widehat{s(D)}\hookrightarrow\widehat{D}$ sends $x^m_i$ to $y^m_i$ for all $i,m$. The inclusion map	$\widehat{t(D)}\hookrightarrow\widehat{D}$ sends $x^{k-1}_i$ to $t(y^k_p)$, where $p=\max\{j:v^k_j(D)=i\}$, if $i\in\im (v^k(D))$. It sends $x^m_i$ to $y^m_i$ whenever $m<k-1$ or $i\notin \im(v^{m+1}(D))$.
	
\end{definition}	

\begin{definition}
	
	Given a simple $k$-string diagram diagram $D$ and a $k$-globular set $X$, we define an \textbf{$X$-labelling} of $D$ to be a map of $k$-globular sets $\widehat{D}\to X$.
	
\end{definition}

This corresponds to the idea of putting labels on the pictorial representation of the string diagram. 

\begin{example}
	
An $X$-labelling of $$\includegraphics[scale=1.5,align=c]{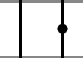}$$ consists of choices of $0$-cells $x, y, z$, $1$-cells $f\colon x\to y$, $g,h\colon y\to z$ and a $2$-cell $\eta\colon g\to h$. We can represent the labeled diagram by inserting the labels into the pictorial representation of the diagram: $$\includegraphics[scale=1.5,align=c]{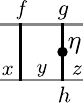}$$

\end{example}

\begin{definition}
	
Let $D$ be a simple $k$-string diagram and $X$ a $k$-globular set. We define source and target maps $$s\colon \gSet_k(\widehat{D},X)\to \gSet_{k-1}(\widehat{s(D)},X)$$ $$t\colon \gSet_k(\widehat{D},X)\to \gSet_{k-1}(\widehat{t(D)},X)$$ by restricting along the inclusion maps $\widehat{s(D)}\hookrightarrow\widehat{D}$ and $\widehat{s(D)}\hookrightarrow\widehat{D}$.	
	
\end{definition}	

Since $\widehat{D}$ is defined as a sequence of pushouts, in order to give a labeling of $D$ it is enough to give labels for each of the cells added via these pushouts. So a labelling $x\colon \widehat{D}\to X$ is equivalent to a choice of $x^m_i\in X_m$ for each $m=0,\cdots, k$ and $i=1,\cdots,\ell_m$ such that \begin{equation*}
s(x^m_i) = \begin{cases}
x^{m-1}_{v^m(i)} &\text{if $i=\min\{j:v^m(j)=v^m(i)\}$}\\
t(x^m_p) &\text{otherwise, with $p=\max\{j<i:v^m(j)=v^m(i)\}$}.
\end{cases}
\end{equation*}  

\begin{example}
	
An $X$-labelling of $$\includegraphics[scale=1.5,align=c]{stringdiag/whisker.pdf}$$ is equivalent to the choice of $x\in X_0$, $f,g\in X_1$ and $\eta\in X_2$ such that $s(f)=x$, $s(g)=t(f)$ and $s(\eta)=g$.	
	
\end{example}

We can specify an $X$-labeling of a simple $k$-string diagram $D$ by assigning labels to an even smaller subset of the cells of $D$, as follows. 

\begin{notation}

An \textbf{$n$-graded set} is a collection $(X(m))_{m=0,\cdots,n}$ of sets. We denote the catgory of graded sets and grade preserving maps by $\Set_{n}$. We sometimes denote $X(k)$ by $X_k$ when there is no risk of confusion with other subscripts and we write $|x|=k$ when $x\in X_k$.
 
\end{notation}

\begin{definition}
	
Given an $n$-globular set $Y$, we define its \textbf{minimal generating subset} to be the $n$-graded subset $\min(Y)\subset Y$ with $\min(Y)_i\subset Y_i$ equal to the complement of the image of the map \[\xymatrixcolsep{4pc}\xymatrix@1{Y_{i+1}\amalg Y_{i+1}\ar[r]^-{s_i\amalg t_i} & Y_i}.\]	
	
\end{definition}	

One can easily check that $\min(Y)$ is the smallest $n$-graded subset of $Y$ such that any $n$-globular subset $Z\subset Y$ containing $\min(Y)$ must be equal to all of $Y$. Moreover, the restriction map $\gSet_n(Y,X)\to\Set_{n}(\min(Y),X)$ is injective.

\begin{remark}

If $D$ is a simple $k$-string diagram, then $$\min(\widehat{D})(m)=\{x^m_i:i\notin\im (v^{m+1}(D))\},$$ where $x^m_i$ are the cells attached in the definition of $\widehat{D}$. 

\end{remark}

\begin{definition}
	
Let $D$ be a simple $k$-string diagram and $X$ a $k$-globular set. A \textbf{minimal $X$-labeling} of $D$ is a map of $n$-graded sets $\min(\widehat{D})\to X$ which is in the image of $\gSet_k(\widehat{D},X)\to\Set_{k}(\min(\widehat{D}),X).$
	
\end{definition}	

The image of the map $\gSet_k(\widehat{D},X)\to\Set_k(\min(\widehat{D}),X)$ is defined by some equations relating the sources and targets of labels of cells in $\min(\widehat{D})$.

\begin{example}
	
An $X$-labelling of $$\includegraphics[scale=1.5,align=c]{stringdiag/whisker.pdf}$$ is determined by the choice of a $1$-cell $f$ and a $2$-cell $\eta$ in $X$, such that $t(f)=s^2(\eta)$. This follows from the fact that the associated globular set is $\C_2\cup_{\C_0}\C_1$. We can depict this labeled diagram by inserting the labels into the pictorial representation of the diagram: $$\includegraphics[scale=1.5,align=c]{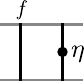}$$
	
\end{example}	

\subsection{The monad $T_n^{D^s}$}

\begin{definition}
	
	We define the functor $T_n^{\sd}\colon \gSet_n\to\gSet_n$ by $$T_n^{\sd}(X)_{k}:=\coprod_{D\in\sd_n(k)}\gSet_n(\widehat{D},X),$$ with source and target maps given by restriction of labellings to $s(D)$ and $t(D)$ respectively. 
	
\end{definition}		

This functor can be equipped with the structure of a monad in much the same way as the monad $T_n$ based on globular pasting diagrams. To define the multiplication map $$\mu^{\sd}\colon T_n^{D^s}T_n^{D^s}\to T_n^{D^s}$$ we need a rule that takes a simple string diagram labelled by simple string diagrams and produces a composite simple string diagram. This can be done by replacing each dot with the diagram that labels it, just as in the case of pasting diagrams. 

For a formal description, let $k\leq n$, $D$ a $k$-diagram, and write $\ell_i=\ell_i(D)$ and $v^i=v^i(D)$. Suppose $D$ is labeled by $m$-diagrams $D^m_i$, for $m=0,\cdots, k$ and $i=1,\cdots,\ell_m$.  Recall the compatibility conditions on the $D^m_i$: \begin{equation*}
s(D^m_i) = \begin{cases}
D^{m-1}_{v^m(i)} &\text{if $i=\min\{j:v^m(j)=v^m(i)\}$}\\
t(D^m_p) &\text{otherwise, with $p=\max\{j<i:v^m(j)=v^m(i)\}$}.
\end{cases}
\end{equation*} Now $t(D^m_p)=s(D^m_p)$ so by induction we always have $s(D^m_i)=D^{m-1}_{v^m(i)}$. Therefore we must have $\ell_j(D^{m-1}_{v^m(i)})=\ell_j(D^m_i)$ and $v^j(D^{m-1}_{v^m(i)})=v^j(D^m_i)$ for $j\leq m-1$. Now we can define the composite diagram to be \[\xymatrix@1{\coprod_{i=1}^{\ell_k}\langle\ell_k(D^k_i)\rangle\ar[r] & \coprod_{i=1}^{\ell_{k-1}}\langle\ell_{k-1}(D^{k-1}_i)\rangle\ar[r] & \cdots \ar[r] & \coprod_{i=1}^{\ell_1}\langle\ell_1(D^1_i)\rangle\ar[r] & \coprod_{i=1}^{\ell_0}\langle\ell_0(D^0_i)\rangle=\langle1\rangle}\] where the maps are \[\xymatrixcolsep{4pc}\xymatrix@1{\langle\ell_m(D^m_i)\rangle\ar[r]^-{v^m(D^m_i)} & \langle\ell_{m-1}(D^m_i)\rangle=\langle\ell_{m-1}(D^{m-1}_{v^m(i)})}.\] The total order on $\coprod_{i=1}^{\ell_m}\langle\ell_m(D^m_i)\rangle$ extends the order on each $\langle\ell_m(D^m_i)\rangle$ by letting every element in the $i$ component be smaller than every element in the $j$ component, whenever $i<j$.

Now suppose the diagrams $D^m_i$ are themselves labeled in some $n$-globular set $X$, in a compatible way, so that we have an element in $T_n^{D^s}T_n^{D^s}(X)$. These labellings correspond to maps $(x^m_i)^p:\langle\ell_p(D^m_i)\rangle\to X_p$ for $p=0,\cdots, m$, satisfying certain compatibility conditions. In particular, we get maps $(x^m_i)^m:\langle\ell_m(D^m_i)\rangle\to X_m$. These induce maps $\coprod_{i=1}^{\ell_m}\langle\ell_m(D^m_i)\rangle\to X_m$, which satisfy the necessary compatibility conditions to give a labelling of the composite diagram, as one can easily check. We thus obtain a map $T_n^{D^s}T_n^{D^s}(X)\to T_n^{D^s}(X)$. One can check that the associativity of this map follows from that of $\amalg$. 

\begin{example}
	
Consider the following $D^s$-labeled simple $3$-string diagram:

$$\includegraphics[scale=1.5,align=c]{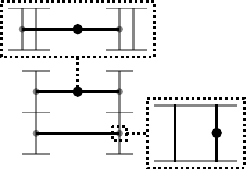}$$ The labels are inside the dashed boxes, connected to the cell which they are labelling. The composite of this diagram is the following: $$\includegraphics[scale=1.5,align=c]{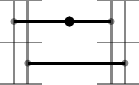}$$	
	
\end{example}	

Finally, we define the unit map $\eta^{\sd}_X\colon X\to T^{D^s}_n(X)$ by sending $x\in X_k$ to the simple $k$-string diagram $I_k:=((1),\cdots,(1))$ with the $k$-cell labelled by $x$. One can easily check that this satisfies the unit axioms, and so we have constructed a monad structure on $T^{D^s}_n$.

\begin{example}
	
Suppose $X$ is an $n$-globular set and $x\in X_4$. Then the unit map $X\to T^{D^s}_nX$ sends $x$ to $$\includegraphics[scale=1.5,align=c]{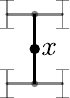}.$$	
	
\end{example}	

\begin{definition}
	
	An \textbf{$n$-sesquicategory} is an algebra over the monad $T_n^{D^s}$. A \textbf{functor} between $n$-sesquicategories is a morphism of algebras over $T_n^{D^s}$. We denote the category $\Alg_{T_n^{D^s}}$ of $T_n^{D^s}$-algebras and $T_n^{D^s}$-algebra morphisms by $\Sesq_n$. 
	
\end{definition} 

Note that what is called a \textbf{sesquicategory} in \cite{street_sesquicat} corresponds in our language to a $2$-sesquicategory. A $1$-sesquicategory is simply a category.

\begin{remark}
	
If we denote by $*$ the $n$-globular set with one element in each dimension, then $T_n^{\sd}(*)=\sd_n$ becomes a globular operad in the sense of \cite{OperadsCats}, when equipped with the map $\pi\colon \sd_n\to\pd_n$ and the product map $T_n(\sd_n)\times_{\pd_n}\sd_n\to\sd_n$ induced by $\mu^{\sd}_{*}$. We have $T_n^{\sd}(X)=T_n(X)\times_{\pd_n}\sd_n$, which means $T_n^{\sd}$ is the monad associated to this globular operad.
	
\end{remark}	

\section{Generators and relations}

An $n$-sesquicategory is an $n$-globular set $\CC$ equipped with an algebraic structure which determines a composite $k$-morphism in $\CC$ for each $\CC$-labelled simple $k$-string diagram. Now we will show that this algebraic structure is determined by certain binary operations of composition and whiskering, together with unit morphisms in each dimension, satisfying certain associativity and unitality equations. 

\subsection{The generators}

We now specify families of simple string diagrams which encode binary operations of composition and whiskering, as well as unary operations encoding units, on any $n$-sesquicategory.

\begin{definition}
	
	A simple $k$-string diagram $D$ is called \textbf{type $i$ binary} when $\ell_i(D)=2$ and $\ell_j(D)=1$ for all $j\neq i$.
	
\end{definition}

There is a unique type $k$ binary simple $k$-string diagram, denoted by $\circ_{k,k}$. We have $\widehat{\circ_{k,k}}=\C_k\cup_{C_{k-1}}\C_k$, so this determines the operation below, usually called \textbf{composition}. $$\circ^{\CC}_{k,k}\colon \CC_k\times_{\CC_{k-1}}\CC_k\to\CC_k$$  

\begin{example}
	
	\begin{center}\begin{tabular}{ccc}$\circ_{1,1}=(2)=\includegraphics[scale=1,align=c]{stringdiag/1d.pdf}$ & ; & $\circ_{2,2}=(1,(1,1))=\includegraphics[scale=1,align=c]{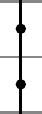}$ .\end{tabular}\end{center}	
	
\end{example}	

For each $i<k$ there are two type $i$ binary simple $k$-string diagrams, which we denote by $\circ_{i,k}$ and $\circ_{k,i}$. They are distinguished by the value of $v^{i+1}$, which is $(1)$ for the first $(2)$ for the second. We have $\widehat{\circ_{i,k}}=\C_i\cup_{C_{i-1}}\C_k$ and $\widehat{\circ_{k,i}}=\C_k\cup_{C_{i-1}}\C_i$, so these determine operations \begin{center}$\circ^{\CC}_{i,k}\colon \CC_i\times_{\CC_{i-1}}\CC_k\to\CC_k$ and $\circ^{\CC}_{k,i}\colon \CC_k\times_{\CC_{i-1}}\CC_i\to\CC_k$ \end{center} which are usually called \textbf{whiskering}. In all cases the pullback is taken by using the source map on the first factor and the target map on the second.

\begin{example} These are the two type $1$ binary simple $2$-string diagrams:
	
	\begin{center}\begin{tabular}{ccc}$\circ_{2,1}=(2,(2))=\includegraphics[scale=1,align=c]{stringdiag/whisker.pdf}$ & ; & $\circ_{1,2}=(2,(1))=\includegraphics[scale=1,align=c]{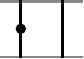}$ .\end{tabular}\end{center}	
	
\end{example}

The simple $k$-string diagram $u_k=(1,(1),\cdots,(1),())$ gives a map $$u_k^{\CC}\colon \CC_{k-1}\to\CC_k.$$ Given a $(k-1)$-morphism $f$ in $\CC$, the $k$-morphism $u_k(f)$ is called the \textbf{identity} on $f$ and also denoted $\id_f$.

\begin{example}
	
	\begin{center}\begin{tabular}{ccc}$u_1=(0)=\includegraphics[scale=1,align=c]{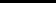}$ & ; & $u_2=(1,())=\includegraphics[scale=1,align=c]{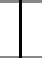}$ .\end{tabular}\end{center}
	
\end{example}

\begin{definition}

We denote by $\OO_n$ the graded subset of $\sd_n$ consisting of the string diagrams $\circ_{i,j}$ and $u_i$, for $i,j=1,\cdots,n$. We refer to the diagrams in $\OO_n$ as \textbf{generating diagrams}.

\end{definition}

\subsection{The relations}\label{relations}

We now describe the relations satisfied by the composition, whiskering and unit operations. These hold in $\sd_n$ and consequently in any $T_n^{\sd}$-algebra.
\begin{definition}
	
	A simple $k$-string diagram $D$ is called \textbf{type $(i,j)$ ternary}, for $i\leq j$, when all the following conditions hold:
	
	\begin{itemize}
		\item if  $i=j$ then $\ell_i(D)=3$ and $\ell_m(D)=1$ for all $m\neq i$;
		\item  if $j\geq i+1$ then $\ell_i(D)=\ell_j(D)=2$ and $\ell_m(D)=1$ for all other $m$;
		\item if $j=i+1$ then $v^j(D)=(1,1)$ or $(2,2)$;

	\end{itemize} 
	
\end{definition}	

As the name suggests, these diagrams determine ternary operations on $\CC$. Moreover, one can construct each of these diagrams by composing binary diagrams in two different ways. This gives equalities between different ways of composing the binary operations, which we interpret as \textbf{associativity} equations.

There is a \textbf{unique type $(k,k)$ ternary simple $k$-string diagram}, namely $$\circ_{k,k,k}:=(1,(1),\cdots,(1),(1,1,1)).$$ This determines the ternary operation $\circ^{\CC}_{k,k,k}\colon \CC_k\times_{\CC_{k-1}}\CC_k\times_{\CC_{k-1}}\CC_k\to\CC_k$, which is equal to both maps $\CC_k\times_{\CC_{k-1}}\CC_k\times_{\CC_{k-1}}\CC_k\to\CC_k$ one can produce by applying $\circ^{\CC}_{k,k}$ twice. This means $\circ_{k,k}^{\CC}$ is associative: $$(x\circ_{k,k}^{\CC} y)\circ_{k,k}^{\CC} z=\circ_{k,k,k}^{\CC}(x,y,z)=x\circ_{k,k}^{\CC} (y\circ_{k,k}^{\CC} z).$$

\begin{example} These are the type $(k,k)$ ternary simple $k$-string diagrams, for $k=1,2$:
	
	\begin{center}\begin{tabular}{ccc}$\circ_{1,1,1}=(3)=\includegraphics[scale=1,align=c]{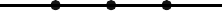}$ & ; & $\circ_{2,2,2}=(1,(1,1,1))=\includegraphics[scale=1,align=c]{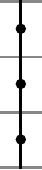}$ .\end{tabular}\end{center} The associativity equation corresponding to $\circ_{1,1,1}$ follows from these identities:
	
	\begin{center}\begin{tabular}{c}$\includegraphics[scale=1]{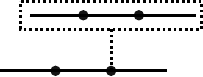}=\includegraphics[scale=1]{stringdiag/111.pdf}=\includegraphics[scale=1]{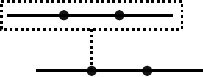}$ \end{tabular}\end{center}
	
\end{example}	

\textbf{For each $i<k$ there are three type $(i,i)$ ternary simple $k$-string diagrams}. They are distinguished only by the value of $v^{i+1}$, which can be $(1)$, $(2)$ or $(3)$, with the corresponding diagram denoted by $\circ_{i,i,k}$, $\circ_{i,k,i}$ and $\circ_{k,i,i}$ respectively. These determine ternary operations $\circ^{\CC}_{i,i,k}$, $\circ^{\CC}_{i,k,i}$ and $\circ^{\CC}_{k,i,i}$. Each of these can be expressed in two ways in terms of the binary operations, and so we get the associativity equations below, relating $\circ^{\CC}_{i,i}$, $\circ^{\CC}_{i,k}$ and $\circ^{\CC}_{k,i}$.

$$(x\circ_{i,i}^{\CC} y)\circ_{i,k}^{\CC} z=\circ_{i,i,k}^{\CC}(x,y,z)=x\circ_{i,k}^{\CC} (y\circ_{i,k}^{\CC} z),$$
$$(x\circ_{i,k}^{\CC} y)\circ_{k,i}^{\CC} z=\circ_{i,k,i}^{\CC}(x,y,z)=x\circ_{i,k}^{\CC} (y\circ_{k,i}^{\CC} z)$$ 
$$(x\circ_{k,i}^{\CC} y)\circ_{k,i}^{\CC} z=\circ_{k,i,i}^{\CC}(x,y,z)=x\circ_{k,i}^{\CC} (y\circ_{i,i}^{\CC} z).$$

\begin{example} These are the three type $(1,1)$ ternary simple $2$-string diagrams:
	
	\begin{center}\begin{tabular}{ccccc}$\circ_{2,1,1}=(3,(3))=\includegraphics[scale=1,align=c]{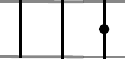}$ & ; & $\circ_{1,2,1}=(3,(2))=\includegraphics[scale=1,align=c]{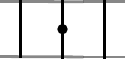}$ & ; & $\circ_{1,1,2}=(3,(1))=\includegraphics[scale=1,align=c]{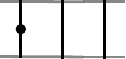}$ .\end{tabular}\end{center} The associativity equation corresponding to $\circ_{2,1,1}$ follows from this identity:
	
	\begin{center}\begin{tabular}{c}$\includegraphics[scale=1]{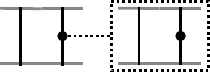}=\includegraphics[scale=1]{stringdiag/211.pdf}=\includegraphics[scale=1]{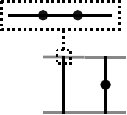}$ \end{tabular}\end{center}
	
\end{example}

\textbf{For each $i<k$ there are two type $(i,k)$ ternary simple $k$-string diagrams}, distinguished by the value of $v^{i+1}$, and denoted by $\circ_{i,k,k}$ and $\circ_{k,k,i}$. These lead to the associativity equations below.
$$x\circ_{i,k}^{\CC} (y\circ_{k,k}^{\CC} z)=\circ_{i,k,k}^{\CC}(x,y,z)=(x\circ_{i,k}^{\CC} y)\circ_{k,k}^{\CC} (x\circ_{i,k}^{\CC} z)$$ $$(x\circ_{k,k}^{\CC} y)\circ_{k,i}^{\CC} z=\circ_{k,k,i}^{\CC}(x,y,z)=(x\circ_{k,i}^{\CC} z)\circ_{k,k}^{\CC} (y\circ_{k,i}^{\CC} z).$$ 

\begin{example} These are the two type $(1,2)$ ternary simple $2$-string diagrams:
	
	\begin{center}\begin{tabular}{ccc}$\circ_{2,2,1}=(2,(2,2))=\includegraphics[scale=1,align=c]{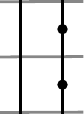}$ & ; & $\circ_{1,2,2}=(2,(1,1))=\includegraphics[scale=1,align=c]{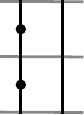}$  .\end{tabular}\end{center} The associativity equation corresponding to $\circ_{2,2,1}$ follows from this identity:
	
	\begin{center}\begin{tabular}{c}$\includegraphics[scale=1,align=c]{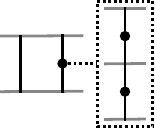}=\includegraphics[scale=1,align=c]{stringdiag/221.pdf}=\includegraphics[scale=1,align=c]{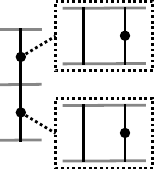}$ \end{tabular}\end{center}
	
\end{example}

\textbf{For each $i<j<k$ there are four type $(i,j)$ ternary simple $k$-string diagrams}, distinguished by the values of $v^{i+1}$ and $v^{j+1}$, and denoted by $\circ_{i,j,k}$, $\circ_{i,k,j}$, $\circ_{j,k,i}$ and $\circ_{k,j,i}$. These lead to the associativity equations below. $$x\circ_{i,k}^{\CC} (y\circ_{j,k}^{\CC} z)=\circ_{i,j,k}^{\CC}(x,y,z)=(x\circ_{i,j}^{\CC} y)\circ_{j,k}^{\CC} (x\circ_{i,k}^{\CC} z),$$  
$$x\circ_{i,k}^{\CC} (y\circ_{k,j}^{\CC} z)=\circ_{i,k,j}^{\CC}(x,y,z)=(x\circ_{i,k}^{\CC} y)\circ_{k,j}^{\CC} (x\circ_{i,j}^{\CC} z),$$  
$$(x\circ_{j,k}^{\CC} y)\circ_{k,i}^{\CC} z=\circ_{j,k,i}^{\CC}(x,y,z)=(x\circ_{j,i}^{\CC} z)\circ_{j,k}^{\CC} (y\circ_{k,i}^{\CC} z)$$ $$(x\circ_{k,j}^{\CC} y)\circ_{k,i}^{\CC} z=\circ_{k,j,i}^{\CC}(x,y,z)=(x\circ_{k,i}^{\CC} z)\circ_{k,j}^{\CC} (y\circ_{j,i}^{\CC} z).$$  

\begin{example} These are the four type $(1,2)$ ternary simple $3$-string diagrams:
	
	\begin{center}\begin{tabular}{ccc}$\circ_{1,2,3}=(2,(1,1),(1))=\includegraphics[scale=1,align=c]{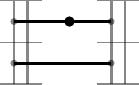}$ & ; & $\circ_{2,1,3}=(2,(2,2),(1))=\includegraphics[scale=1,align=c]{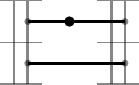}$  \\ \\ $\circ_{3,1,2}=(2,(1,1),(2))=\includegraphics[scale=1,align=c]{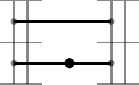}$  & ; & $\circ_{3,2,1}=(2,(2,2),(2))=\includegraphics[scale=1,align=c]{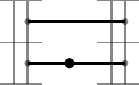}$  .\end{tabular}\end{center} The associativity equation corresponding to $\circ_{1,2,3}$ follows from this identity:
	
	\begin{center}\begin{tabular}{c}$\includegraphics[scale=1]{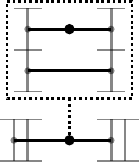}=\includegraphics[scale=1]{stringdiag/123.pdf}=\includegraphics[scale=1]{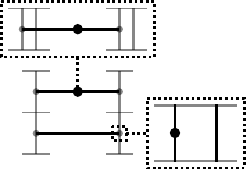}$ \end{tabular}\end{center}

\end{example}

Finally, we have the following \textbf{unitality} relations, which can easily be checked by composing diagrams. \begin{center}\begin{tabular}{ll}$(i\leq k):$& $u_i^{\CC}(x)\circ_{i,k}^{\CC}y=y\text{ and }x\circ_{k,i}^{\CC}u_i^{\CC}(y)=x$ \\ \\ $(i<k):$ & $u_k^{\CC}(x)\circ_{k,i}^{\CC}y=u_k(x\circ_{k-1,i}^{\CC}y)\text{ and }x\circ_{i,k}^{\CC}u_k^{\CC}(y)=u_k(x\circ_{i,k-1}^{\CC}y)$ \end{tabular}\end{center} 

\begin{definition}
	
We denote by $\EE_n$ the set of relations described above, involving generators of dimension $\leq n$.
	
\end{definition}	

\section{The monad $T_n^{\sd}$ is generated by $(\OO_n,\EE_n)$}

In this section, we give an explicit construction of a monad $\widetilde{\Tree}_n^{\OO,\EE}$ generated by $(\OO_n,\EE_n)$ and we show that this monad is isomorphic to $T_n^{\sd}$. We also include a remark relating $\widetilde{\Tree}_n^{\OO,\EE}$ to the globular operad generated by a presentation in the sense of \cite{PresGlobOp}.

\subsection{The monad $\widetilde{\Tree}_n^{\OO,\EE}$}

We now define a monad $\widetilde{\Tree}_n^{\OO,\EE}$, which is generated by the generators and relations we have just described.

\begin{notation}
	
By a \textbf{rooted tree} $(T,r)$ we mean a connected graph with no cycles with a chosen distinguished vertex $r$, which we call its \textbf{root}. We denote by $V(T)$ and $E(T)$ the sets of vertices and edges of $T$, respectively. We orient the edges of $T$ towards $r$. We denote the source and target maps by $s,t:E(T)\to V(T)$. A \textbf{leaf} of $T$ is a vertex with no incoming edges. We denote the set of leaves of $T$ by $L(T)$. An \textbf{internal vertex} of $T$ is a vertex which is not a leaf. We denote the set of internal vertices of $T$ by $I(T)$. Given $v\in V(T)$ we denote by $E^-(v)$ the set of \textbf{incoming edges} at $v$. The \textbf{height} of $T$ is the maximum length of a directed path in $T$.
	
\end{notation}	

\begin{definition}
	
Let $X$ be an $n$-graded set. A $k$-dimensional $(\OO_n,X)$-labelled tree is a rooted tree $(T,r)$, together with

\begin{enumerate}
	
\item a map $\lambda\colon I(T)\to\OO_n$;
\item a map $\lambda\colon L(T)\to X$;
\item a bijection $\phi_v\colon \min(\widehat{\lambda(v)})\to E^-(v)$ for each $v\in I(T)$;		
\end{enumerate}	 	

such that

\begin{enumerate}
	
 \item $|\lambda(r)|=k;$
 \item for each $v\in I(T)$, the map $\xymatrixcolsep{1.5pc}\xymatrix{\min(\widehat{\lambda(v)})\ar[r]^-{\phi_v} & E^-(v)\ar[r]^-{s} & V(T)\ar[r]^-{\lambda} & X\coprod\OO_n}$ is graded.
 
\end{enumerate} 

We denote the set of $k$-dimensional $(\OO_n,X)$-labeled trees by $\Tree_n^{\OO}(X)_{k}$ or $\Tree_n^{\OO}(X)_k$.
	
\end{definition}

So internal vertices are labelled by diagrams in $\OO_n$ and leaves are labelled by cells in $X$. The root label must have dimension $k$. At each internal vertex $v$ with labelling diagram $D$, we have an incoming edge for each cell in $\min(\widehat{D})$ and moreover the label associated to the source of each of these edges must have dimension equal to that of the corresponding cell. The idea is that this will give a minimal labeling of $D$, once certain compatibility conditions are satisfied.

\begin{notation}
	
If $x\in\Tree_n^{\OO}(X)_i$ and $y\in\Tree_n^{\OO}(X)_{j}$, then we denote by \[\xymatrixcolsep{1pc}\xymatrixrowsep{1pc}\vcenter{\vbox{\xymatrix{ & \circ_{i,j} &  \\ x\ar[ru] & & y\ar[lu]\\}}}\] the $\max\{i,j\}$-dimensional $(\OO_n,X)$-labelled tree obtained by connecting the roots of $x$ and $y$ to a new common root $r$, with $\lambda(r)=\circ_{i,j}$. Similarly, if $z\in\Tree_n^{\OO}(X)(i-1)$, then $z\to u_{i}$ is an $i$-dimensional $(\OO_n,X)$-labelled tree. We call this process \textbf{grafting}.
	
\end{notation}

\begin{definition}
	
	An $n$-\textbf{preglobular set} is an $n$-graded set $X$ together with maps $s,t\colon X_k\to X_{k-1}$ for each $k\geq 1$.	
	
\end{definition}
 
\begin{definition}
	
Let $X$ be an $n$-preglobular set. We define source and target maps $$s,t\colon \Tree_n^{\OO}(X)_{k}\to \Tree_n^{\OO}(X)(k-1).$$ For trees of height zero, these are the source and target maps of $X$. For trees of nonzero height, we use the following inductive formulas for $s$, where $i<k$ and $x$ and $y$ have appropriate dimensions in each case. The map $t$ is defined by the same formulas, replacing every instance of $s$ with $t$.

\begin{center}\begin{tabular}{l}

$s(\xymatrix{x\ar[r] & u_k})=x;$ 

\\ \\

$s\left(\vcenter{\vbox{\xymatrix{x\ar[r] & \circ_{k,k} & y\ar[l]}}}\right)=s(y);$ 

\\ \\

 $s\left(\vcenter{\vbox{\xymatrix{x\ar[r] & \circ_{i,k} & y\ar[l]}}}\right)=\vcenter{\vbox{\xymatrix{x\ar[r] & \circ_{i,k-1} & s(y)\ar[l]\\}}};$ 
 
 \\ \\

 $s\left(\vcenter{\vbox{\xymatrix{x\ar[r] & \circ_{k,i} & y\ar[l]}}}\right)=\vcenter{\vbox{\xymatrix{s(x)\ar[r] & \circ_{k-1,i} & y\ar[l]}}};$ 
 
 \end{tabular}\end{center}

\end{definition}	

Note that, even if we assume that $X$ is an $n$-globular set, these source and target maps don't necessarily satisfy the globularity conditions, so $\Tree_n^{\OO}(X)$ is in general only an $n$-preglobular set. 

\begin{notation}
	
We denote by $\tau_{\leq h}\Tree_n^{\OO}(X)\subset\Tree_n^{\OO}(X)$ the $n$-preglobular subset consisting of trees of height $\leq h$.	
	
\end{notation}

\begin{definition}
	
Let $X$ be an $n$-preglobular set. A \textbf{globular relation} on $X$ is an equivalence relation $\sim$ on each $X_k$ such that

\begin{enumerate}
	
\item if $x\sim\tilde{x}$ then $s(x)\sim s(\tilde{x})$ and $t(x)\sim t(\tilde{x})$;

\item $ss(x)\sim st(x)$ and $ts(x)\sim tt(x)$ for all x.	
	
\end{enumerate}		
	
\end{definition}

This means that the quotient $(X\quot\sim)$ is an $n$-globular set.

\begin{definition} Consider an $n$-preglobular set $X$ with globular relation $\sim$. We define, by induction on $h$, an $n$-preglobular subset $\tau_{\leq h}\Tree_n^{\OO,\EE}(X,\sim)\subset\tau_{\leq h}\Tree_n^{\OO}(X)$ with a relation $\stackrel{\epsilon}{\sim}_h$. Elements in $\tau_{\leq h}\Tree_n^{\OO,\EE}(X,\sim)$ are called \textbf{$\stackrel{\epsilon}{\sim}_{h-1}$-compatible}. The inductive definition is presented below.

\end{definition} When $h=0$, we let $\tau_{\leq 0}\Tree_n^{\OO,\EE}(X,\sim):=\tau_{\leq 0}\Tree_n^{\OO}(X)=X$ and the relation $\stackrel{\epsilon}{\sim}_0$ is $\sim$. 

Now consider $h\geq 1$. Any $x\in\tau_{\leq h}\Tree_n^{\OO}(X)$ of height zero is $\stackrel{\epsilon}{\sim}_{h-1}$-compatible. Let $x\in\tau_{\leq h-1}\Tree_n^{\OO}(X)_i$, $y\in\tau_{\leq h-1}\Tree_n^{\OO}(X)_{j}$ and $m=\min\{i,j\}$. Then $$\xymatrixcolsep{1pc}\xymatrixrowsep{1pc}\xymatrix{ & \circ_{i,j} &  \\ x\ar[ru] & & y\ar[lu]}$$  is $\stackrel{\epsilon}{\sim}_{h-1}$-compatible if and only if $x,y$ are $\stackrel{\epsilon}{\sim}_{h-2}$-compatible and $s^{i-m+1}(x)\stackrel{\epsilon}{\sim}_{h-1}t^{j-m+1}(y)$. Moreover, $x\to u_{i+1}$ is $\stackrel{\epsilon}{\sim}_{h-1}$-compatible if and only if $x$ is $\stackrel{\epsilon}{\sim}_{h-2}$-compatible. Now we must define the globular relation $\stackrel{\epsilon}{\sim}_{h}$ on $\tau_{\leq h}\Tree_n^{\OO,\EE}(X,\sim)$.

If $x,y\in\tau_{\leq h}\Tree_n^{\OO,\EE}(X,\sim)$ have height zero and $x\stackrel{\epsilon}{\sim}_{0}y$, then $x\stackrel{\epsilon}{\sim}_{h}y$. 

Let $i\leq k$, $x\in \tau_{\leq h-2}\Tree_n^{\OO,\EE}(X,\sim)_{i-1}$, $y\in \tau_{\leq h-1}\Tree_n^{\OO,\EE}(X,\sim)_k$. If $x\stackrel{\epsilon}{\sim}_{h-1} t^{k-i+1}(y)$, then \begin{center}\begin{tabular}{lccc} $(\lambda_{i,k}):$ & $\vcenter{\vbox{\xymatrixcolsep{1pc}\xymatrixrowsep{1pc}\xymatrix{ & \circ_{i,k} &  \\ u_i\ar[ru] & & y\ar[lu]\\ x\ar[u] & & }}}$ & $\stackrel{\epsilon}{\sim}_{h}$ & $y$.\end{tabular}\end{center}

Let $i\leq k$, $x\in \tau_{\leq h-1}\Tree_n^{\OO,\EE}(X,\sim)_k$, $y\in \tau_{\leq h-2}\Tree_n^{\OO,\EE}(X,\sim)_{i-1}$. If $s^{k-i+1}(x)\stackrel{\epsilon}{\sim}_{h-1}y$, then \begin{center}\begin{tabular}{lccc} $(\rho_{k,i}):$ & $\vcenter{\vbox{\xymatrixcolsep{1pc}\xymatrixrowsep{1pc}\xymatrix{ & \circ_{k,i} &  \\ x\ar[ru] & & u_i\ar[lu]\\  & & y\ar[u]}}}$ & $\stackrel{\epsilon}{\sim}_{h}$ & $x$.\end{tabular}\end{center}

Let $i<k$, $x\in \tau_{\leq h-2}\Tree_n^{\OO,\EE}(X,\sim)_i$, $y\in \tau_{\leq h-2}\Tree_n^{\OO,\EE}(X,\sim)_{k-1}$. If $s(x)\stackrel{\epsilon}{\sim}_{h-2}t^{k-i}(y)$, then \begin{center}\begin{tabular}{lccc} $(\rho_{i,k}):$ & $\vcenter{\vbox{\xymatrixcolsep{1pc}\xymatrixrowsep{1pc}\xymatrix{ & \circ_{i,k} &  \\ x\ar[ru] & & u_k\ar[lu]\\ & & y\ar[u]}}}$ & $\stackrel{\epsilon}{\sim}_{h}$ & $\vcenter{\vbox{\xymatrixcolsep{1pc}\xymatrixrowsep{1pc}\xymatrix{ & u_k &  \\  & \circ_{i,k-1}\ar[u] & \\ x\ar[ru] & & y\ar[lu]}}}$.\end{tabular}\end{center}

Let $i<k$, $x\in \tau_{\leq h-2}\Tree_n^{\OO,\EE}(X,\sim)_{k-1}$, $y\in \tau_{\leq h-2}\Tree_n^{\OO,\EE}(X,\sim)_i$. If $s^{k-i}(x)\stackrel{\epsilon}{\sim}_{h-2}t(y)$, then \begin{center}\begin{tabular}{lccc} $(\lambda_{k,i}):$ & $\vcenter{\vbox{\xymatrixcolsep{1pc}\xymatrixrowsep{1pc}\xymatrix{ & \circ_{k,i} &  \\ u_k\ar[ru] & & y\ar[lu]\\ x\ar[u] & & }}}$ & $\stackrel{\epsilon}{\sim}_{h}$ & $\vcenter{\vbox{\xymatrixcolsep{1pc}\xymatrixrowsep{1pc}\xymatrix{ & u_k &  \\  & \circ_{k-1,i}\ar[u] & \\ x\ar[ru] & & y\ar[lu]}}}$.\end{tabular}\end{center}

Let $k\geq 1$ and $x,y,z\in \tau_{\leq h-2}\Tree_n^{\OO,\EE}(X,\sim)_k$. If $s(x)\stackrel{\epsilon}{\sim}_{h-2}t(y)$ and $s(y)\stackrel{\epsilon}{\sim}_{h-2}t(z)$, then \begin{center}\begin{tabular}{lccc} $(\circ_{k,k,k}):$ & $\vcenter{\vbox{\xymatrixcolsep{1pc}\xymatrixrowsep{1pc}\xymatrix{ & & \circ_{k,k} &  \\ & \circ_{k,k}\ar[ru] & & z\ar[lu] \\ x\ar[ru] & & y\ar[lu] & }}}$ & $\stackrel{\epsilon}{\sim}_{h}$ & $\vcenter{\vbox{\xymatrixcolsep{1pc}\xymatrixrowsep{1pc}\xymatrix{ & \circ_{k,k} & &  \\ x\ar[ru] & & \circ_{k,k}\ar[lu] & \\ & y\ar[ru] & & z\ar[lu]}}}$.\end{tabular}\end{center}

Let $i<k$, $x,y\in \tau_{\leq h-2}\Tree_n^{\OO,\EE}(X,\sim)_i$, $z\in \tau_{\leq h-2}\Tree_n^{\OO,\EE}(X,\sim)_k$. If $s(x)\stackrel{\epsilon}{\sim}_{h-2}t(y)$ and $s(y)\stackrel{\epsilon}{\sim}_{h-2}t^{k-i+1}(z)$, then \begin{center}\begin{tabular}{lccc} $(\circ_{i,i,k}):$ & $\vcenter{\vbox{\xymatrixcolsep{1pc}\xymatrixrowsep{1pc}\xymatrix{ & & \circ_{i,k} &  \\ & \circ_{i,i}\ar[ru] & & z\ar[lu] \\ x\ar[ru] & & y\ar[lu] & }}}$ & $\stackrel{\epsilon}{\sim}_{h}$ & $\vcenter{\vbox{\xymatrixcolsep{1pc}\xymatrixrowsep{1pc}\xymatrix{ & \circ_{i,k} & &  \\ x\ar[ru] & & \circ_{i,k}\ar[lu] & \\ & y\ar[ru] & & z\ar[lu]}}}$.\end{tabular}\end{center}

Let $i<k$, $x,z\in \tau_{\leq h-2}\Tree_n^{\OO,\EE}(X,\sim)_i$, $y\in \tau_{\leq h-2}\Tree_n^{\OO,\EE}(X,\sim)_k$. If $s(x)\stackrel{\epsilon}{\sim}_{h-2}t^{k-i+1}(y)$ and $s^{k-i+1}(y)\stackrel{\epsilon}{\sim}_{h-2}t(z)$, then \begin{center}\begin{tabular}{lccc} $(\circ_{i,k,i}):$ & $\vcenter{\vbox{\xymatrixcolsep{1pc}\xymatrixrowsep{1pc}\xymatrix{ & & \circ_{k,i} &  \\ & \circ_{i,k}\ar[ru] & & z\ar[lu] \\ x\ar[ru] & & y\ar[lu] & }}}$ & $\stackrel{\epsilon}{\sim}_{h}$ & $\vcenter{\vbox{\xymatrixcolsep{1pc}\xymatrixrowsep{1pc}\xymatrix{ & \circ_{i,k} & &  \\ x\ar[ru] & & \circ_{k,i}\ar[lu] & \\ & y\ar[ru] & & z\ar[lu]}}}$.\end{tabular}\end{center}

Let $i<k$, $x\in \tau_{\leq h-2}\Tree_n^{\OO,\EE}(X,\sim)_k$, $y,z\in \tau_{\leq h-2}\Tree_n^{\OO,\EE}(X,\sim)_i$. If $s^{k-i+1}(x)\stackrel{\epsilon}{\sim}_{h-2}t(y)$ and $s(y)\stackrel{\epsilon}{\sim}_{h-2}t(z)$, then \begin{center}\begin{tabular}{lccc} $(\circ_{k,i,i}):$ & $\vcenter{\vbox{\xymatrixcolsep{1pc}\xymatrixrowsep{1pc}\xymatrix{ & & \circ_{k,i} &  \\ & \circ_{k,i}\ar[ru] & & z\ar[lu] \\ x\ar[ru] & & y\ar[lu] & }}}$ & $\stackrel{\epsilon}{\sim}_{h}$ & $\vcenter{\vbox{\xymatrixcolsep{1pc}\xymatrixrowsep{1pc}\xymatrix{ & \circ_{k,i} & &  \\ x\ar[ru] & & \circ_{i,i}\ar[lu] & \\ & y\ar[ru] & & z\ar[lu]}}}$.\end{tabular}\end{center}

Let $i<k$, $x\in \tau_{\leq h-2}\Tree_n^{\OO,\EE}(X,\sim)_i$, $y,z\in \tau_{\leq h-2}\Tree_n^{\OO,\EE}(X,\sim)_k$. If $s(x)\stackrel{\epsilon}{\sim}_{h-2}t^{k-i+1}(y)$ and $s(y)\stackrel{\epsilon}{\sim}_{h-2}t(z)$, then \begin{center}\begin{tabular}{lccc} $(\circ_{i,k,k}):$ & $\vcenter{\vbox{\xymatrixcolsep{1pc}\xymatrixrowsep{1pc}\xymatrix{ & \circ_{i,k} & &  \\ x\ar[ru] & & \circ_{k,k}\ar[lu] & \\ & y\ar[ru] & & z\ar[lu]}}}$ & $\stackrel{\epsilon}{\sim}_{h}$ & $\vcenter{\vbox{\xymatrixcolsep{1pc}\xymatrixrowsep{1pc}\xymatrix{ & & \circ_{k,k} & & \\ & \circ_{i,k}\ar[ru] & & \circ_{i,k}\ar[lu] & \\ x\ar[ru] & & y\ar[lu]\text{  }x\ar[ru] & & z\ar[lu]}}}$.\end{tabular}\end{center}

Let $i<k$, $x,y\in \tau_{\leq h-2}\Tree_n^{\OO,\EE}(X,\sim)_k$, $z\in \tau_{\leq h-2}\Tree_n^{\OO,\EE}(X,\sim)_i$. If $s(x)\stackrel{\epsilon}{\sim}_{h-2}t(y)$ and $s^{k-i+1}(y)\stackrel{\epsilon}{\sim}_{h-2}t(z)$, then \begin{center}\begin{tabular}{lccc} $(\circ_{k,k,i}):$ & $\vcenter{\vbox{\xymatrixcolsep{1pc}\xymatrixrowsep{1pc}\xymatrix{ & & \circ_{k,i} &  \\ & \circ_{k,k}\ar[ru] & & z\ar[lu] \\ x\ar[ru] & & y\ar[lu] & }}}$ & $\stackrel{\epsilon}{\sim}_{h}$ & $\vcenter{\vbox{\xymatrixcolsep{1pc}\xymatrixrowsep{1pc}\xymatrix{ & & \circ_{k,k} & & \\ & \circ_{k,i}\ar[ru] & & \circ_{k,i}\ar[lu] & \\ x\ar[ru] & & z\ar[lu]\text{  }y\ar[ru] & & z\ar[lu]}}}$.\end{tabular}\end{center}

Let $i<j<k$ and take $x\in \tau_{\leq h-2}\Tree_n^{\OO,\EE}(X,\sim)_i$, $y\in \tau_{\leq h-2}\Tree_n^{\OO,\EE}(X,\sim)_j$ and $z\in \tau_{\leq h-2}\Tree_n^{\OO,\EE}(X,\sim)_k$. If $s(x)\stackrel{\epsilon}{\sim}_{h-2}t^{j-i+1}(y)$ and $s(y)\stackrel{\epsilon}{\sim}_{h-2}t^{k-j+1}(z)$, then \begin{center}\begin{tabular}{lccc} $(\circ_{i,j,k}):$ & $\vcenter{\vbox{\xymatrixcolsep{1pc}\xymatrixrowsep{1pc}\xymatrix{ & \circ_{i,k} & &  \\ x\ar[ru] & & \circ_{j,k}\ar[lu] & \\ & y\ar[ru] & & z\ar[lu]}}}$ & $\stackrel{\epsilon}{\sim}_{h}$ & $\vcenter{\vbox{\xymatrixcolsep{1pc}\xymatrixrowsep{1pc}\xymatrix{ & & \circ_{j,k} & & \\ & \circ_{i,j}\ar[ru] & & \circ_{i,k}\ar[lu] & \\ x\ar[ru] & & y\ar[lu]\text{  }x\ar[ru] & & z\ar[lu]}}}$.\end{tabular}\end{center}

Let $i<j<k$ and take $x\in \tau_{\leq h-2}\Tree_n^{\OO,\EE}(X,\sim)_i$, $y\in \tau_{\leq h-2}\Tree_n^{\OO,\EE}(X,\sim)_k$ and $z\in \tau_{\leq h-2}\Tree_n^{\OO,\EE}(X,\sim)_j$. If $s(x)\stackrel{\epsilon}{\sim}_{h-2}t^{k-i+1}(y)$ and $s^{k-j+1}(y)\stackrel{\epsilon}{\sim}_{h-2}t(z)$, then \begin{center}\begin{tabular}{lccc} $(\circ_{i,k,j}):$ & $\vcenter{\vbox{\xymatrixcolsep{1pc}\xymatrixrowsep{1pc}\xymatrix{ & \circ_{i,k} & &  \\ x\ar[ru] & & \circ_{k,j}\ar[lu] & \\ & y\ar[ru] & & z\ar[lu]}}}$ & $\stackrel{\epsilon}{\sim}_{h}$ & $\vcenter{\vbox{\xymatrixcolsep{1pc}\xymatrixrowsep{1pc}\xymatrix{ & & \circ_{k,j} & & \\ & \circ_{i,k}\ar[ru] & & \circ_{i,j}\ar[lu] & \\ x\ar[ru] & & y\ar[lu]\text{  }x\ar[ru] & & z\ar[lu]}}}$.\end{tabular}\end{center}

Let $i<j<k$ and take $x\in \tau_{\leq h-2}\Tree_n^{\OO,\EE}(X,\sim)_j$, $y\in \tau_{\leq h-2}\Tree_n^{\OO,\EE}(X,\sim)_k$ and $z\in \tau_{\leq h-2}\Tree_n^{\OO,\EE}(X,\sim)_i$. If $s(x)\stackrel{\epsilon}{\sim}_{h-2}t^{k-j+1}(y)$ and $s^{k-i+1}(y)\stackrel{\epsilon}{\sim}_{h-2}t(z)$, then \begin{center}\begin{tabular}{lccc} $(\circ_{j,k,i}):$ & $\vcenter{\vbox{\xymatrixcolsep{1pc}\xymatrixrowsep{1pc}\xymatrix{ & & \circ_{k,i} &  \\ & \circ_{j,k}\ar[ru] & & z\ar[lu] \\ x\ar[ru] & & y\ar[lu] & }}}$ & $\stackrel{\epsilon}{\sim}_{h}$ & $\vcenter{\vbox{\xymatrixcolsep{1pc}\xymatrixrowsep{1pc}\xymatrix{ & & \circ_{j,k} & & \\ & \circ_{j,i}\ar[ru] & & \circ_{k,i}\ar[lu] & \\ x\ar[ru] & & z\ar[lu]\text{  }y\ar[ru] & & z\ar[lu]}}}$.\end{tabular}\end{center}

Let $i<j<k$ and take $x\in \tau_{\leq h-2}\Tree_n^{\OO,\EE}(X,\sim)_k$, $y\in \tau_{\leq h-2}\Tree_n^{\OO,\EE}(X,\sim)_j$ and $z\in \tau_{\leq h-2}\Tree_n^{\OO,\EE}(X,\sim)_i$. If $s^{k-j+1}(x)\stackrel{\epsilon}{\sim}_{h-2}t(y)$ and $s^{j-i+1}(y)\stackrel{\epsilon}{\sim}_{h-2}t(z)$, then \begin{center}\begin{tabular}{lccc} $(\circ_{k,j,i}):$ & $\vcenter{\vbox{\xymatrixcolsep{1pc}\xymatrixrowsep{1pc}\xymatrix{ & & \circ_{k,i} &  \\ & \circ_{k,j}\ar[ru] & & z\ar[lu] \\ x\ar[ru] & & y\ar[lu] & }}}$ & $\stackrel{\epsilon}{\sim}_{h}$ & $\vcenter{\vbox{\xymatrixcolsep{1pc}\xymatrixrowsep{1pc}\xymatrix{ & & \circ_{k,j} & & \\ & \circ_{k,i}\ar[ru] & & \circ_{j,i}\ar[lu] & \\ x\ar[ru] & & z\ar[lu]\text{  }y\ar[ru] & & z\ar[lu]}}}$.\end{tabular}\end{center}

Let $x,\tilde{x}\in \tau_{\leq h-1}\Tree_n^{\OO,\EE}(X,\sim)_{k-1}$. If $x\stackrel{\epsilon}{\sim}_{h-1} \tilde{x}$, then \begin{center}\begin{tabular}{lccc} $(u_k):$ & $\vcenter{\vbox{\xymatrixcolsep{1pc}\xymatrixrowsep{1pc}\xymatrix{u_k \\ x\ar[u]}}}$ & $\stackrel{\epsilon}{\sim}_{h}$ &  $\vcenter{\vbox{\xymatrixcolsep{1pc}\xymatrixrowsep{1pc}\xymatrix{u_k \\ \tilde{x}\ar[u]}}}$.\end{tabular}\end{center}

Let $x,\tilde{x}\in \tau_{\leq h-1}\Tree_n^{\OO,\EE}(X,\sim)_i$, $y,\tilde{y}\in\tau_{\leq h-1}\Tree_n^{\OO,\EE}(X,\sim)_j$ and $m=\min\{i,j\}$. If $x\stackrel{\epsilon}{\sim}_{h-1} \tilde{x}$, $y\stackrel{\epsilon}{\sim}_{h-1} \tilde{y}$, $s^{i-m+1}(x)\stackrel{\epsilon}{\sim}_{h-1} t^{j-m+1}(y)$ and $s^{i-m+1}(\tilde{x})\stackrel{\epsilon}{\sim}_{h-1} t^{j-m+1}(\tilde{y})$ then \begin{center}\begin{tabular}{lccc} $(\circ_{i,j}):$ & $\vcenter{\vbox{\xymatrixcolsep{1pc}\xymatrixrowsep{1pc}\xymatrix{ & \circ_{i,j} &  \\ x\ar[ru] & & y\ar[lu]}}}$ & $\stackrel{\epsilon}{\sim}_{h}$ & $\vcenter{\vbox{\xymatrixcolsep{1pc}\xymatrixrowsep{1pc}\xymatrix{ & \circ_{i,j} &  \\ \tilde{x}\ar[ru] & & \tilde{y}\ar[lu]}}}$.\end{tabular}\end{center}  

\begin{remark}\label{conditions}

Note that for each of the equations defining $\stackrel{\epsilon}{\sim}_h$ there are some conditions on subtrees $x,y,z$ appearing on both sides of the equation. In each case, these conditions imply that both sides of the equation are are $\stackrel{\epsilon}{\sim}_{h-1}$-compatible, so these equations do define a relation on $\tau_{\leq h}\Tree_n^{\OO,\EE}(X)$.

Moreover, if the left hand side of an equation is $\stackrel{\epsilon}{\sim}_{h-1}$-compatible, then the conditions are satisfied. This means that when we have a $\stackrel{\epsilon}{\sim}_{h-1}$-compatible tree we can apply these equations from left to right without checking any conditions.

In all cases except $(\circ_{i,k,k}), (\circ_{k,k,i}), (\circ_{i,j,k}), (\circ_{i,k,j}), (\circ_{j,k,i})$ and $(\circ_{k,j,i})$, if the right hand side of the equation is $\stackrel{\epsilon}{\sim}_{h-1}$-compatible, then the conditions are satisfied, so we can also apply the equations from right to left without checking any conditions. 
\end{remark}	

\begin{remark}\label{hh+1}
	
Clearly, if $x$ is $\stackrel{\epsilon}{\sim}_{h-1}$-compatible, then it is also $\stackrel{\epsilon}{\sim}_{h}$-compatible. Moreover, if $x\stackrel{\epsilon}{\sim}_h y$ then $x\stackrel{\epsilon}{\sim}_{h+1}y$.

\end{remark}	

\begin{lemma}\label{hpreglob}

Let $(X,\sim)$ be an $n$-preglobular set with globular relation and let $x,y\in\tau_{\leq h}\Tree_n^{\OO}(X)$.

\begin{enumerate}
 \item If $x$ is $\stackrel{\epsilon}{\sim}_{h-1}$-compatible, then $s(x), t(x)$ are $\stackrel{\epsilon}{\sim}_{h-1}$-compatible.

\item If $x\stackrel{\epsilon}{\sim}_{h} y$ then $s(x)\stackrel{\epsilon}{\sim}_{h}s(y)$ and $t(x)\stackrel{\epsilon}{\sim}_{h}t(y)$.

\item If $x$ is $\stackrel{\epsilon}{\sim}_{h-1}$-compatible, then $ss(x)\stackrel{\epsilon}{\sim}_{h}st(x)$ and $ts(x)\stackrel{\epsilon}{\sim}_{h}tt(x)$.

\end{enumerate}

This implies that $(\tau_{\leq h}\Tree_n^{\OO,\EE}(X,\sim),\stackrel{\epsilon}{\sim}_{h})$ is an $n$-preglobular set with globular relation.
 
\end{lemma}

\begin{proof}
 
The proof is by induction on $h$. The case $h=0$ follows directly from the fact that $(X,\sim)$ is an $n$-preglobular set with globular relation. Now let $h\geq 1$. Starting with \emph{1.}, there is one case for each possible labeling of the root and in each such case one only needs to perform a simple computation using the formulas defining $s$ and $t$. For example, suppose $x=\vcenter{\vbox{\xymatrixcolsep{1pc}\xymatrixrowsep{1pc}\xymatrix{y\ar[r] & \circ_{i,k} & z\ar[l]}}},$ so that $s(x)=\vcenter{\vbox{\xymatrixcolsep{1pc}\xymatrixrowsep{1pc}\xymatrix{ y\ar[r] & \circ_{i,k-1} & s(z)\ar[l]}}}.$ Since $x$ is $\stackrel{\epsilon}{\sim}_{h-1}$-compatible, we have $y,z$ $\stackrel{\epsilon}{\sim}_{h-2}$-compatible. By induction this implies that is $s(z)$ is $\stackrel{\epsilon}{\sim}_{h-2}$-compatible. So all we have left to show is that $s(y)\stackrel{\epsilon}{\sim}_{h-1}t^{k-1-i+1}s(z)$. But $t^{k-1-i+1}s(z)\stackrel{\epsilon}{\sim}_{h-1}t^{k-i+1}(z)$ by induction (using \emph{2} and \emph{3}), and we have $s(y)\stackrel{\epsilon}{\sim}_{h-1}t^{k-i+1}(z)$ because $x$ is $\stackrel{\epsilon}{\sim}_{h-1}$-compatible. 

For \emph{2.}, we need to do a simple computation for each of the equations defining $\stackrel{\epsilon}{\sim}_h$. For example, consider the equation \begin{center}\begin{tabular}{lccc} $(\circ_{k,k,i}):$ & $\vcenter{\vbox{\xymatrixcolsep{1pc}\xymatrixrowsep{1pc}\xymatrix{ & & \circ_{k,i} &  \\ & \circ_{k,k}\ar[ru] & & z\ar[lu] \\ x\ar[ru] & & y\ar[lu] & }}}$ & $\stackrel{\epsilon}{\sim}_h$ & $\vcenter{\vbox{\xymatrixcolsep{1pc}\xymatrixrowsep{1pc}\xymatrix{ & & \circ_{k,k} & & \\ & \circ_{k,i}\ar[ru] & & \circ_{k,i}\ar[lu] & \\ x\ar[ru] & & z\ar[lu]\text{  }y\ar[ru] & & z\ar[lu]}}}$.\end{tabular}\end{center} Both sides have source $\vcenter{\vbox{\xymatrixcolsep{1pc}\xymatrixrowsep{1pc}\xymatrix{ s(y) \ar[r] & \circ_{k-1,i} & z\ar[l]}}}.$ 	

For \emph{3.}, we need to do a computation for each possible labeling of the root. For example, we have $$ss\left(\vcenter{\vbox{\xymatrixcolsep{1pc}\xymatrixrowsep{1pc}\xymatrix{x\ar[r] & \circ_{i,i+1} & y\ar[l] }}}\right)=s\left(\vcenter{\vbox{\xymatrixcolsep{1pc}\xymatrixrowsep{1pc}\xymatrix{x\ar[r] & \circ_{i,i} & s(y)\ar[l]}}}\right)=s^2(y)$$ $$st\left(\vcenter{\vbox{\xymatrixcolsep{1pc}\xymatrixrowsep{1pc}\xymatrix{x\ar[r] & \circ_{i,i+1} & y\ar[l]}}}\right)=s\left(\vcenter{\vbox{\xymatrixcolsep{1pc}\xymatrixrowsep{1pc}\xymatrix{x\ar[r] & \circ_{i,i} & t(y)\ar[l]}}}\right)=st(y)$$ and we have $ss(y)\stackrel{\epsilon}{\sim}_{h-1}st(y)$ by induction.\end{proof}

\begin{definition}
	
Let $(X,\sim)$ be an $n$-preglobular set  with a globular relation. An element $x\in\Tree_n^{\OO}(X)$ is called $\stackrel{\epsilon}{\sim}$-compatible when it is $\stackrel{\epsilon}{\sim}_{h}$-compatible for some $h$. We denote the set of $k$-dimensional $\stackrel{\epsilon}{\sim}$-compatible $(\OO_n,X)$-labeled trees by $\Tree_n^{\OO,\EE}(X,\sim)_{k}.$ We define a relation $\stackrel{\epsilon}{\sim}$ on $\Tree_n^{\OO,\EE}(X,\sim)$ by letting $x\stackrel{\epsilon}{\sim}y$ if and only if $x\stackrel{\epsilon}{\sim}_{h} y$ for some $h$.

\end{definition}

\begin{proposition}
	
Let $(X,\sim)$ be an $n$-preglobular set with a globular relation. Then $(\Tree_n^{\OO,\EE}(X,\sim),\stackrel{\epsilon}{\sim})$ is an $n$-preglobular set with globular relation.
	
\end{proposition}

\begin{proof}

This follows easily from Lemma \ref{hpreglob} and Remark \ref{hh+1}. \end{proof}

\begin{notation}
	
	When $X$ is an $n$-globular set, we write $\Tree_n^{\OO,\EE}(X):=\Tree_n^{\OO,\EE}(X,=).$ This preglobular set is equipped with the globular relation $\stackrel{\epsilon}{=}$. 	
	
\end{notation}	

\begin{definition}
	
Let $(X,\sim)$ be an $n$-preglobular set with globular relation. We write $$\widetilde{\Tree}_n^{\OO,\EE}(X,\sim):=\Tree_n^{\OO,\EE}(X,\sim)\quot\stackrel{\epsilon}{\sim}.$$ When $X$ is an $n$-globular set, we write $\widetilde{\Tree}_n^{\OO,\EE}(X):=\Tree_n^{\OO,\EE}(X)\quot\stackrel{\epsilon}{=}.$
	
\end{definition}

\begin{notation}
	
Let $X$ be an $n$-globular set. We write $$\Tree_n^{\OO,\EE}(\Tree_n^{\OO,\EE}(X)):=\Tree_n^{\OO,\EE}(\Tree_n^{\OO,\EE}(X),\stackrel{\epsilon}{=})$$ and we denote by $\stackrel{\epsilon^2}{=}$ the globular relation on this preglobular set which extends $\stackrel{\epsilon}{=}$. 	
	
\end{notation}	

\begin{definition}

Let $X$ be an $n$-globular set. We define a map $$\mu^{\OO}_X\colon \Tree_n^{\OO}(\Tree_n^{\OO}(X))\to\Tree_n^{\OO}(X)$$ by induction on height. An $x\in\Tree_n^{\OO}(\Tree_n^{\OO}(X))$ of height $0$ is just an element in $\Tree_n^{\OO}(X)$. Given $x$ of height $\geq 1$, let $o\in\OO$ be its root label and let $x_i\in\Tree_n^{\OO}(\Tree_n^{\OO}(X))$ be the subtrees with root at distance $1$ from the root of $x$. Then we obtain $\mu^{\OO}_X(x)$ by grafting $\mu^{\OO}_X(x_i)$ onto $o$. 

\end{definition}

\begin{lemma}\label{mu}
 
The map $\mu^{\OO}_X$ defined above has the following properties.

\begin{enumerate}
 \item If $x$ is $\stackrel{\epsilon^2}{=}$-compatible, then $\mu^{\OO}_X(x)$ is $\stackrel{\epsilon}{=}$-compatible;
 \item if $x,y$ are $\stackrel{\epsilon^2}{=}$-compatible and $x\stackrel{\epsilon^2}{=}y$, then $\mu^{\OO}_X(x)\stackrel{\epsilon}{=}\mu^{\OO}_X(y)$;
 \item $s(\mu^{\OO}_X(x))=\mu^{\OO}_X(s(x))$ and $t(\mu^{\OO}_X(x))=\mu^{\OO}_X(t(x))$.
\end{enumerate}

\end{lemma}

\begin{proof}
 
This is easy to prove by induction.\end{proof}

\begin{example}\label{grafting}
	
We construct $\tau\in \Tree_n^{\OO}(\Tree_n^{\OO,\EE}(X))$ such that $\mu^{\OO}_X(\tau)$ is $\stackrel{\epsilon}{=}$-compatible, while $\tau$ itself is not $\stackrel{\epsilon^2}{=}$-compatible. In other words, we have $\mu^{\OO}_X(\tau)\in\Tree_n^{\OO,\EE}(X)$, but $\tau\notin\Tree_n^{\OO,\EE}(\Tree_n^{\OO,\EE}(X)).$	

Let $A,B,C,D\in\Tree_n^{\OO}(X)$ be defined by $A=x\in X(1)$, $D=r\in X(1)$, \begin{center}\begin{tabular}{ccc}$B=\vcenter{\vbox{\xymatrixcolsep{1pc}\xymatrixrowsep{1pc}\xymatrix{ & \circ_{2,1} &  \\ z\ar[ru] & & y\ar[lu]}}}$ & and &  $C=\vcenter{\vbox{\xymatrixcolsep{1pc}\xymatrixrowsep{1pc}\xymatrix{ & \circ_{2,1} &  \\ q\ar[ru] & & p\ar[lu]}}}$, \end{tabular}\end{center} with $z,q\in X(2)$ and $y,p\in X(1)$. We suppose $s^2(z)=t(y)$ and $s^2(q)=t(p)$, so $B,C$ are $\stackrel{\epsilon}{=}$-compatible.

Now let $\tau$ the $(\OO_n,\Tree_n^{\OO}(X))$-labeled tree $$\vcenter{\vbox{\xymatrixcolsep{1pc}\xymatrixrowsep{1pc}\xymatrix{ & & \circ_{2,2} & & \\ & \circ_{1,2}\ar[ru] & & \circ_{2,1}\ar[lu] & \\ D\ar[ru] & & C\ar[lu]\text{  }B\ar[ru] & & A\ar[lu]}}},$$ whose image in $\Tree_n^{\OO}(X)$ under $\mu^{\OO}_X$ is $$\vcenter{\vbox{\xymatrixcolsep{1pc}\xymatrixrowsep{1pc}\xymatrix{ & & \circ_{2,2} & & \\ & \circ_{1,2}\ar[ru] & & \circ_{2,1}\ar[lu] & \\ r\ar[ru] & \circ_{2,1}\ar[u] & & \circ_{2,1}\ar[u] & x\ar[lu] \\ q\ar[ru] & p\ar[u] & & z\ar[u] & y\ar[lu] }}}.$$ One checks easily that this is $\stackrel{\epsilon}{=}$-compatible if and only if $s(r)=t^2(q)$, $s(y)=t(x)$, $p=x$, $s(q)=y$ and $r=t(z)$.

Now for $\tau$ to be $\stackrel{\epsilon^2}{=}$-compatible, we would need to have $$s\left(\vcenter{\vbox{\xymatrixcolsep{1pc}\xymatrixrowsep{1pc}\xymatrix{ & \circ_{1,2} &  \\ D\ar[ru] & & C\ar[lu]}}}\right)\stackrel{\epsilon^2}{=} t\left(\vcenter{\vbox{\xymatrixcolsep{1pc}\xymatrixrowsep{1pc}\xymatrix{ & \circ_{2,1} &  \\ B\ar[ru] & & A\ar[lu]}}}\right)$$ which is equivalent to $$\vcenter{\vbox{\xymatrixcolsep{1pc}\xymatrixrowsep{1pc}\xymatrix{ & \circ_{1,1} &  \\ D\ar[ru] & & s(C)\ar[lu]}}}\stackrel{\epsilon^2}{=} \vcenter{\vbox{\xymatrixcolsep{1pc}\xymatrixrowsep{1pc}\xymatrix{ & \circ_{1,1} &  \\ t(B)\ar[ru] & & A\ar[lu]}}}$$ which is equivalent to $D\stackrel{\epsilon}{=} t(B)$ and $s(C)\stackrel{\epsilon}{=} A$, which does not hold.	
\end{example}	

We will now equip $\widetilde{\Tree}_n^{\OO,\EE}$ with the structure of a monad on $\gSet_n$.	By Lemma \ref{mu}, the map $\mu^{\OO}_X$ defined above induces a map $$\mu^{\OO}_X\colon \widetilde{\Tree}_n^{\OO,\EE}(\widetilde{\Tree}_n^{\OO,\EE}(X))\to\widetilde{\Tree}_n^{\OO,\EE}(X).$$ We define the map $$\eta^{\OO}_X\colon X\to\widetilde{\Tree}_n^{\OO,\EE}(X)$$ by sending $x\in X$ to the tree with one vertex, labeled by $X$.

\begin{lemma}

The triple $(\widetilde{\Tree}_n^{\OO,\EE},\mu^{\OO},\eta^{\OO})$ is a monad.

\end{lemma}	

\begin{proof}
	
It's easy to show that $\mu^{\OO}_X$ is associative, and unital with respect to $\eta$.	\end{proof}	

\subsection{The map $\tilde{\varphi}\colon \widetilde{\Tree}_n^{\OO,\EE}\to T^{\sd}_n$}

Now we define a map of monads $\tilde{\varphi}\colon \widetilde{\Tree}_n^{\OO,\EE}\to T^{\sd}_n$. This will be obtained by passing to the quotient of the $\stackrel{\epsilon}{=}$-invariant map of $n$-preglobular sets $\varphi_X\colon \Tree_n^{\OO,\EE}(X)\to T_n^{\sd}(X)$ which we define  by induction on the height of the tree, as follows.

\begin{notation}
 
If $X$ is an $n$-globular set and $x\in X_i$, $y\in X_k$, we denote by $\xymatrixcolsep{1.5pc}\xymatrix{x\ar[r] & \circ_{i,k} & y\ar[l]}$, $\xymatrixcolsep{1.5pc}\xymatrix{x\ar[r] & u_{i+1}}$ and $\xymatrixcolsep{1.5pc}\xymatrix{x\ar[r] & I_{i}}$ the obvious elements of $T_n^{\sd}(X)$. We are already using the same notation for trees, but these will be distinguished from context.
 
\end{notation}

\begin{definition}

 We define a map $$(\varphi_X)_0\colon \tau_{\leq 0}\Tree_n^{\OO,\EE}(X)\to T_n^{\sd}(X)$$ by sending $x\in X(k)$ to $x\to I_k$.
 
  \end{definition}
  
This map is $\stackrel{\epsilon}{=}_0$-invariant, since $\stackrel{\epsilon}{=}_0$ is just $=$. It sends $s(x)$ to $s(x)\to I_{k-1}=s(x\to I_{k})$ and $t(x)$ to $t(x)\to I_{k-1}=t(x\to I_k)$, so it is a map of $n$-preglobular sets
 
 \begin{definition}
 
 An $x\in \tau_{\leq 1}\Tree_n^{\OO,\EE}(X)$ of height $1$ is an element $d\in\OO_n(k)$ with an $X$-labeling, which is already an element in $T_n^{\sd}(X)$, so we can extend $(\varphi_X)_0$ to a map $$(\varphi_X)_1\colon \tau_{\leq 1}\Tree_n^{\OO,\EE}(X)\to T_n^{\sd}(X).$$
 
  \end{definition}
  
The relation $\stackrel{\epsilon}{=}_1$ is again just $=$, so $(\varphi_X)_1$ is $\stackrel{\epsilon}{=}_1$-invariant. One also checks easily that it is a map of $n$-preglobular sets.
 
 \begin{definition}
 
 Now let $h>1$ and suppose we have already defined a $\stackrel{\epsilon}{=}_{h-1}$-invariant map $(\varphi_X)_{h-1}\colon \tau_{\leq h-1} \Tree_n^{\OO,\EE}(X)\to T_n^{\sd}(X).$ We define $$(\varphi_X)_{h}\colon \tau_{\leq h} \Tree_n^{\OO,\EE}(X)\to T_n^{\sd}(X)$$ extending $(\varphi_X)_{h-1}$ as follows. Let $x\in\tau_{\leq h} \Tree_n^{\OO,\EE}(X)$ have height $h$, let $(T,r)$ be the underlying rooted tree and let $\lambda(r)=d\in\OO_n(k)$. Apply $(\varphi_X)_{h-1}$ to each of the subtrees whose root is a vertex at distance $1$ from $r$. Since $x$ is $\stackrel{\epsilon}{=}_{h-1}$-compatible and  $(\varphi_X)_{h-1}$ is $\stackrel{\epsilon}{=}_{h-1}$-invariant, we get a $T_n^{\sd}(X)$-labeling of $d$. Now we can use the multiplication map $\mu^{\sd}_X\colon T_n^{\sd}(T_n^{\sd}(X))\to T_n^{\sd}(X)$ to obtain an element in $T_n^{\sd}(X)$. 
 
 \end{definition}
 
 One can easily check that $(\varphi_X)_h$ is a map of $n$-preglobular sets.

 \begin{lemma}\label{phinv}
 	
 Suppose $(\varphi_X)_{h-1}$ is $\stackrel{\epsilon}{=}_{h-1}$-invariant and $(\varphi_X)_{h-2}$ is $\stackrel{\epsilon}{=}_{h-2}$-invariant. Then $(\varphi_X)_h$ is $\stackrel{\epsilon}{=}_{h}$-invariant.	
 	
 \end{lemma}	
 
 \begin{proof}
 	
 We need to check that $(\varphi_X)_h$ assigns the same value to boths sides of each equation defining $\stackrel{\epsilon}{=}_{h}$. For the equations $(\circ_{i,j})$ and $(u_k)$ this follows simply from the fact that $(\varphi_X)_{h-1}$ is $\stackrel{\epsilon}{=}_{h-1}$-invariant. For the rest we have to do a computation. We do it only in the example of $(\circ_{i,k,k})$ as the others are similar. So let $i<k$ and $x\in \tau_{\leq h-2}\Tree_n^{\OO,\EE}(X)_i$, $y,z\in \tau_{\leq h-2}\Tree_n^{\OO,\EE}(X)_{k}$. Suppose $s(x)\stackrel{\epsilon}{=}_{h-2}t^{k-i+1}(y)$ and $s(y)\stackrel{\epsilon}{=}_{h-2}t(z)$ and recall the equation \begin{center}\begin{tabular}{lccc} $(\circ_{i,k,k}):$ & $\vcenter{\vbox{\xymatrixcolsep{1pc}\xymatrixrowsep{1pc}\xymatrix{ & \circ_{i,k} & &  \\ x\ar[ru] & & \circ_{k,k}\ar[lu] & \\ & y\ar[ru] & & z\ar[lu]}}}$ & $\stackrel{\epsilon}{=}_{h}$ & $\vcenter{\vbox{\xymatrixcolsep{1pc}\xymatrixrowsep{1pc}\xymatrix{ & & \circ_{k,k} & & \\ & \circ_{i,k}\ar[ru] & & \circ_{i,k}\ar[lu] & \\ x\ar[ru] & & y\ar[lu]\text{  }x\ar[ru] & & z\ar[lu]}}}$.\end{tabular}\end{center} We denote the left and right hand side by $L,R\in\tau_{\leq h}\Tree_n^{\OO,\EE}(X)$. Let $\tilde{x},\tilde{y},\tilde{z}\in T_n^{\sd}(X)$ be the result of applying $(\varphi_X)_{h-2}$ to $x,y,z$. Since $(\varphi_X)_{h-2}$ is $\stackrel{\epsilon}{=}_{h-2}$-invariant, we have $s(\tilde{x})=t^{k-i+1}(\tilde{y})$ and $s(\tilde{y})=t(\tilde{z})$. This implies that $\tilde{L},\tilde{R}$ defined by \begin{center}\begin{tabular}{ccc}$\tilde{L}:=\vcenter{\vbox{\xymatrixcolsep{1pc}\xymatrixrowsep{1pc}\xymatrix{ & \circ_{i,k} & &  \\ I_i\ar[ru] & & \circ_{k,k}\ar[lu] & \\ \tilde{x}\ar[u] & \tilde{y}\ar[ru] & & \tilde{z}\ar[lu]}}}$ & and & $\tilde{R}:=\vcenter{\vbox{\xymatrixcolsep{1pc}\xymatrixrowsep{1pc}\xymatrix{ & & \circ_{k,k} & & \\ & \circ_{i,k}\ar[ru] & & \circ_{i,k}\ar[lu] & \\ \tilde{x}\ar[ru] & & \tilde{y}\ar[lu]\text{  }\tilde{x}\ar[ru] & & \tilde{z}\ar[lu]}}}$ \end{tabular}\end{center} are elements in $(T_n^{\sd})^3(X)$. Since $(\varphi_X)_h$ is computed by using $\mu^{\sd}$ to compose, from bottom to top, we know that $(\varphi_X)_h(L)$ and $(\varphi_X)_h(R)$ are equal to the result of applying the composite map \[\xymatrix@1{T_n^{\sd}((T_n^{\sd})^2(X))\ar[rr]^-{T_n^{\sd}(\mu^{\sd}_X)} && T_n^{\sd}(T_n^{\sd}(X))\ar[r]^-{\mu_X^{\sd}} & T_n^{\sd}(X)}\] to $\tilde{L}$ and $\tilde{R}$, respectively. Since $\mu^{\sd}$ is associative, this is equal to the composite \[\xymatrix@1{(T_n^{\sd})^2(T_n^{\sd}(X))\ar[rr]^-{\mu^{\sd}_{T_n^{\sd}(X)}} && T_n^{\sd}(T_n^{\sd}(X))\ar[r]^-{\mu_X^{\sd}} & T_n^{\sd}(X)}.\] But the first map in this latter composite sends both $\tilde{L}$ and $\tilde{R}$ to \[\vcenter{\vbox{\xymatrixcolsep{1pc}\xymatrixrowsep{1pc}\xymatrix{ & \circ_{i,k,k} & \\ \tilde{x}\ar[ru] & \tilde{y}\ar[u] & \tilde{z}\ar[lu]}}},\] so $(\varphi_X)_h(L)=(\varphi_X)_h(R)$.\end{proof}	
 
\begin{definition}
	
We define the map of $n$-preglobular sets $$\varphi_X\colon \Tree_n^{\OO,\EE}(X)\to T_n^{\sd}(X)$$ by letting $\varphi_X(x)=(\varphi_X)_h(x)$ where $h$ is the height of $x$.	
	
\end{definition}	

\begin{remark}
	
	Unraveling the inductive definition, we see that $\varphi_X(x)$ can be computed by collapsing the internal edges of $x$ from the leaves up to the root, using $\mu^{\sd}_X$ to compose the diagrams which label the corresponding vertices. 
	
\end{remark}

\begin{lemma}
	
The map $\varphi_X$ defined above is $\stackrel{\epsilon}{=}$-invariant.	
	
\end{lemma}	 

\begin{proof}Follows immediately from Lemma \ref{phinv}.\end{proof}

\begin{definition}
	
We define $$\tilde{\varphi}\colon \widetilde{\Tree}_n^{\OO,\EE}\to T^{\sd}_n$$ to be the map obtained from $\varphi$ by passing to the quotient.	
	
\end{definition}	
 
 Now we want to show that $\tilde{\varphi}$ is a map of monads. 
 
\begin{lemma}
	
	The map $\tilde{\varphi}\colon \widetilde{\Tree}_n^{\OO,\EE}\to T^{\sd}_n$ preserves the monad units.
	
\end{lemma} 

\begin{proof}
	
We must show that $\tilde{\varphi}_X\circ\eta^{\OO}_X=\eta^{\sd}_X$. This is clear from the definition of $(\varphi_X)_0$.\end{proof}	

We will need the following notation in the proof that $\tilde{\varphi}$ preserves monad multiplication.

\begin{notation}
	
We denote by $\tau_h\Tree_n^{\OO,\EE}(X,\sim)\subset\tau_{\leq h}\Tree_n^{\OO,\EE}(X,\sim)$ the subset of trees with height $h$.	
	
\end{notation}

\begin{remark}

When $X$ is an $n$-globular set we have an inclusion map $$\tau_1\Tree_n^{\OO,\EE}(X)\hookrightarrow T_n^{\sd}(X)$$ given by interpreting a $(\OO_n,X)$-labeled $\stackrel{\epsilon}{=}_0$-compatible tree of height $1$ as an $X$-labeling of the diagram assigned to the root. Moreover, there is a map $$\tau_h\Tree_n^{\OO,\EE}(X)\to\tau_1\Tree_n^{\OO,\EE}(\tau_{\leq h-1}\Tree_n^{\OO,\EE}(X),\stackrel{\epsilon}{=}),$$ which sends an $(\OO_n,X)$-labelled $\stackrel{\epsilon}{=}_{h-1}$-compatible tree $(T,r)$ of height $h$ to the subtree consisting only of the vertices at distance $\leq 1$ from $r$. The label for the root remains unchanged and each vertex at distance $1$ from the root is labelled by the subtree of $T$ whose root is that vertex. This map has the property that the composite \[\xymatrix{\tau_h\Tree_n^{\OO,\EE}(X)\ar[r] & \tau_1\Tree_n^{\OO,\EE}(\tau_{\leq h-1}\Tree_n^{\OO,\EE}(X),\stackrel{\epsilon}{=}) \ar[r] & \tau_{\leq h}\Tree_n^{\OO,\EE}(X)}\] with the grafting map is equal to the inclusion $\tau_{h}\Tree_n^{\OO,\EE}(X)\hookrightarrow \tau_{\leq h}\Tree_n^{\OO,\EE}(X)$.

\end{remark}

\begin{lemma}
	
The map $\tilde{\varphi}\colon \widetilde{\Tree}_n^{\OO,\EE}\to T^{\sd}_n$ defined above preserves monad multiplication.	
	
\end{lemma}	

\begin{proof}
	
Given an $n$-globular set $X$, we must show that the composites \[\xymatrixcolsep{1.5pc}\xymatrix@1{\Tree_n^{\OO,\EE}(\Tree_n^{\OO,\EE}(X),\stackrel{\epsilon}{=})\ar[rrr]^-{\Tree_n^{\OO,\EE}(\varphi_X)} &&& \Tree_n^{\OO,\EE}(T_n^{\sd}(X))\ar[rr]^-{\varphi_{T_n^{\sd}(X)}} && T_n^{\sd}(T_n^{\sd}(X))\ar[r]^-{\mu^{\sd}_X} & T_n^{\sd}(X) }\] and \[\xymatrix@1{\Tree_n^{\OO,\EE}(\Tree_n^{\OO,\EE}(X),\stackrel{\epsilon}{=})\ar[r]^-{\mu^{\OO}_X} & \Tree_n^{\OO,\EE}(X)\ar[r]^-{\varphi_X} & T_n^{\sd}(X)}\] are equal. We will do this by showing, by induction on $h$, that the composites \[\xymatrixcolsep{1.5pc}\xymatrix@1{\tau_h\Tree_n^{\OO,\EE}(\Tree_n^{\OO,\EE}(X),\stackrel{\epsilon}{=})\ar[rrr]^-{\tau_h\Tree_n^{\OO,\EE}(\varphi_X)} &&& \tau_h\Tree_n^{\OO,\EE}(T_n^{\sd}(X))\ar[rr]^-{\varphi_{T_n^{\sd}(X)}} && T_n^{\sd}(T_n^{\sd}(X))\ar[r]^-{\mu^{\sd}_X} & T_n^{\sd}(X) }\] and \[\xymatrix@1{\tau_h\Tree_n^{\OO,\EE}(\Tree_n^{\OO,\EE}(X),\stackrel{\epsilon}{=})\ar[r]^-{\mu^{\OO}_X} & \Tree_n^{\OO,\EE}(X)\ar[r]^-{\varphi_X} & T_n^{\sd}(X)}\] are equal for all $h$. We denote these composites by $\alpha_h$ and $\beta_h$, respectively. For $h=0$, there is nothing to prove. When $h=1$, both composites are simply given by collapsing all edges from the leaves up to the root, using $\mu^{\sd}$ to compose the corresponding diagrams. 

Now suppose $h\geq2$. We show $\alpha_h$ and $\beta_h$ are equal to the top and bottom composites, respectively, in the following diagram: 

\[\xymatrix@1{ & & (T_n^{\sd})^2(X)\ar[rd]^{\mu^{\sd}_X} & \\ \tau_h\Tree_n^{\OO,\EE}(\Tree_n^{\OO,\EE}(X),\stackrel{\epsilon}{=}) \ar[r]^-{\psi_h} & (T_n^{\sd})^3(X)\ar[ru]^{\mu^{\sd}_{T_n^{\sd}(X)}}\ar[rd]_{T_n^{\sd}(\mu_X^{\sd})}  & &T_n^{\sd}(X). \\ & & (T_n^{\sd})^2(X)\ar[ru]_{\mu^{\sd}_X} &}\] Note that these two composites are equal, by associativity of $\mu^{\sd}$. The map $\psi_h$ is the following composite: 

\[\xymatrixcolsep{2pc}\xymatrix@1{\tau_h\Tree_n^{\OO,\EE}(\Tree_n^{\OO,\EE}(X),\stackrel{\epsilon}{=}) \ar[rrr]^-{\tau_h\Tree_n^{\OO,\EE}(\varphi_X)} &&& \tau_h\Tree_n^{\OO,\EE}(T_n^{\sd}(X))\ar[llld]  \\ \tau_1\Tree_n^{\OO,\EE}(\tau_{\leq h-1}\Tree_n^{\OO,\EE}(T_n^{\sd}(X)),\stackrel{\epsilon}{=}) \ar[rrr]_-{\tau_1\Tree_n^{\OO,\EE}(\varphi_{T_n^{\sd}(X)})} &&& \tau_1\Tree_n^{\OO,\EE}((T_n^{\sd})^2(X))\text{ } \ar@{^{(}->}[r] & (T_n^{\sd})^3(X) .}\] Given $x\in\tau_h\Tree_n^{\OO,\EE}(\Tree_n^{\OO,\EE}(X),\stackrel{\epsilon}{=})$ with underlying tree $(T,r)$, the map $\psi_h$ collapses the edges in the trees which label the leaves of $T$ and then the edges connecting internal vertices of $T$, excluding the ones which connect to the root $r$. In this way one can see directly that the top composite is equal to $\alpha_h$. 

Now consider the map $T_n^{\sd}(\mu_X^{\sd})\circ\psi_h$, i.e., the composition the first two maps on the bottom of the diagram. This consists in applying the map $\alpha_g$ to each of the subtrees of $T$ rooted at a vertex at distance $1$ from $r$, where $g$ is the height of the subtree in question, so $g<h$. By induction on $h$, this is equal to applying $\beta_g$ to each such subtree. Composing this with $\mu_X^{\sd}$ we therefore obtain $\beta_h$, as desired.\end{proof}	

So we have proved the following result.

\begin{proposition}

The map $\tilde{\varphi}\colon \widetilde{\Tree}_n^{\OO,\EE}\to T^{\sd}_n$ is a map of monads.

\end{proposition}

\begin{remark}
	
Let $*$ denote the $n$-globular set with one cell in each dimension. The composite \[\xymatrixcolsep{2pc}\xymatrix@1{\widetilde{\Tree}_n^{\OO,\EE}(*)\ar[r]^{\tilde{\varphi}} & T_n^{\sd}(*)\ar[r]^{\pi} & T_n(*)}\] makes $\widetilde{\Tree}_n^{\OO,\EE}(*)$ an $n$-globular operad, in the sense of \cite{OperadsCats}. We have $$\widetilde{\Tree}_n^{\OO,\EE}(X)=T_n(X)\times_{T_n(*)}\widetilde{\Tree}_n^{\OO,\EE}(*)$$ so $\widetilde{\Tree}_n^{\OO,\EE}$ is the monad associated to this operad. Given an $n$-globular operad $\GG_n$ and a graded map $\theta\colon \OO_n\to\GG_n$ such that the relations $\EE_n$ are satisfied in $\GG_n$, there is a unique map of $n$-globular operads $\widetilde{\Tree}_n^{\OO,\EE}(*)\to\GG_n$ extending $\theta$. We can construct this map in exactly the same way we constructed $\tilde{\varphi}$. This means that $\widetilde{\Tree}_n^{\OO,\EE}(*)$ is the operad generated by the generators $\OO_n$ and the relations $\EE_n$, in the sense of \cite{PresGlobOp}	
	
\end{remark}	

This Remark already makes clear that algebras over $\widetilde{\Tree}_n^{\OO,\EE}$ are exactly the objects described in Theorem \ref{genrel}.

\subsection{$\tilde{\varphi}$ is injective}

Now we prove that the map $\tilde{\varphi}\colon \widetilde{\Tree}_n^{\OO,\EE}\to T^{\sd}_n$ is injective. 

\begin{definition}
	
An $(\OO_n,X)$-labeled tree is \textbf{nondegenerate} if it doesn't contain any instances of $u_i$.

\end{definition}

\begin{lemma}\label{nondegenerate1}
	
Any $x\in\Tree_n^{\OO,\EE}(X)$ is $\stackrel{\epsilon}{=}$-equivalent to a tree of the form $y\to u_{k-i}\to \cdots \to u_{k}$ where $y$ is nondegenerate of dimension $k-i-1$.  	
	
\end{lemma}	

\begin{proof}
	
One can use $(\lambda)$ and $(\rho)$ relations to eliminate some $u_j$ generators and then to push others up towards the root.\end{proof}	

\begin{lemma}\label{nondegenerate2}
	
Let $y_1,y_2\in\Tree_n^{\OO,\EE}(X)$ be nondegenerate and suppose $$\varphi_X(y_1\to u_{k-i_1}\to \cdots \to u_{k})=\varphi_X(y_2\to u_{k-i_2}\to \cdots \to u_{k}).$$ Then $i_1=i_2$ and $\varphi_X(y_1)=\varphi_X(y_2)$.
	
\end{lemma}	

\begin{proof}
	
Write $D=\varphi_X(y_1\to u_{k-i_1}\to \cdots \to u_{k})=\varphi_X(y_2\to u_{k-i_2}\to \cdots \to u_{k})$. Then $i_1=i_2=\max\{i:\ell_{k-i}(D)=0\}$ and $\varphi_X(y_1)=\varphi_X(y_2)=(v^1(D),\cdots,v^{k-i-1}(D))$.\end{proof}	

So we reduced the problem to showing that $\varphi$ is injective on nondegenerate trees.

\begin{definition}
	
A nondegenerate $(\OO_n,X)$-labeled tree is $m$-\textbf{ordered} if for each internal edge $v_1\to v_2$ we have $m(\lambda(v_1))\leq m(\lambda(v_2))$.
	
\end{definition}

\begin{lemma}\label{ordered}
	
Let $x\in\Tree_n^{\OO,\EE}(X)$ and suppose $x$ is nondegenerate. Then $x$ is equivalent to an $m$-ordered tree.	
	
\end{lemma}	

\begin{proof}
	
One can move to an $m$-ordered tree by applying relations $(\circ_{i,k,k})$ and $(\circ_{k,k,i})$, for $i<k$; and $(\circ_{i,j,k}),(\circ_{i,k,j}),(\circ_{j,k,i})$ and $(\circ_{k,j,i})$, for $i<j<k$. Note that in order to move to an $m$-ordered tree we only need to apply these equations from left hand side to right hand side, which means that we don't need to check any conditions, as explained in Remark \ref{conditions}.\end{proof}	

So we reduced the problem to showing that $\varphi$ is injective on $m$-ordered trees. We will do this by induction on the following invariant.

\begin{definition}
	
	Let $x=(T,r,\lambda,\phi)\in\Tree_n^{\OO}(X)$ be nondegenerate. Define $$H(x):=|\{i:\exists_{v\in I(T)}~m(\lambda(v))=i\}|.$$
	
\end{definition}	

We will also need the following definitions.

\begin{definition}

A nondegenerate $x\in\Tree_n^{\OO}(X)$ is called \textbf{nontrivial} if its height is $\geq 1$.

\end{definition}

\begin{definition}
	
Let $x=(T,r,\lambda,\phi)\in\Tree_n^{\OO}(X)$ be nontrivial. We define \begin{center}$m(x)=\min\{m(\lambda(v)):v\in I(T)\}$ and $M(x)=\max\{M(\lambda(v)):v\in I(T)\}.$\end{center}
	
\end{definition}	

\begin{definition}
	
	For any $D\in T_n^{\sd}(X)$, we define $H(D):=|\{i:\ell_i(D)>1\}|$.	
	
\end{definition}

\begin{definition}

A diagram  $D\in T_n^{\sd}(X)$ is \textbf{nondegenerate} if $\ell_i(D)>0$ for all $i$.	

A nondegenerate diagram $D\in T_n^{\sd}(X)$ is \textbf{nontrivial} if $\ell_i(D)>1$ for some $i$.
	
\end{definition}	

\begin{definition}
	
	Given $D\in T_n^{\sd}(X)$ nontrivial, we define \begin{center}$m(D)=\min\{j:\ell_j(D)>1\}$ and $M(D)=\max\{j:\ell_j(D)>1\}.$\end{center}
	
\end{definition}

\begin{remark}
	
	If $D=\circ_{i,j}$, we have $m(D)=M(D)=\min\{i,j\}$.	
	
\end{remark}

\begin{lemma}\label{mMH}
	
Let $x=(T,r,\lambda,\phi)\in\Tree_n^{\OO,\EE}(X)$ be nondegenerate. Then we have $H(\varphi_X(x))=H(x)$. If $x$ has height $\geq 1$, we have $m(\varphi_X(x))=m(x)$, $M(\varphi_X(x))=M(x)$. \end{lemma}	

\begin{proof}
	
One checks, by induction on the height of $x$, that $$\{i:\ell_i(\varphi_X(x))>1\}=\{i:\exists_{v\in I(T)}~m(\lambda(v))=i\}.$$\end{proof}	

\begin{definition}
	
A nontrivial diagram $D\in T_n^{\sd}(X)$ is a \textbf{block} if $m(D)=M(D).$	
	
\end{definition}	

\begin{notation}
	
Let $1\leq i < k$, $\ell>1$ and $1\leq p \leq \ell$. We denote by $B^k(i,\ell,p)$ the simple $k$-string diagram $$B^k(i,\ell,p):=((1),\cdots,(1),(1,\cdots,1),(p),(1),\cdots,(1)),$$ i.e. we have $\ell_i(B^k(i,\ell,p))=\ell$, $\ell_j(B^k(i,\ell,p))=1$ for $j\neq i$ and $v^{i+1}(B^k(i,\ell,p))=(p)$. 

For $k\geq 1$ and $\ell>1$, we denote by $B^{k}(k,\ell)$ the simple $k$-string diagram $$B^{k}(k,\ell):=((1),\cdots,(1),(1,\cdots,1)),$$ i.e.  $\ell_j(B^{k}(k,\ell))=1$ for $j<k$ and $\ell_k(B^{k}(k,\ell))=\ell$.  		

These diagrams are blocks, and every block is of one of these two types.
	
\end{notation}	

\begin{remark}	
	
	The associated globular set of $B^k(i,\ell,p)$ is $$\widehat{B^k}(i,\ell,p)=C_i\cup_{C_{i-1}} \cdots\cup_{C_{i-1}} C_i\cup_{C_{i-1}}C_k\cup_{C_{i-1}}C_i\cup_{C_{i-1}}\cdots\cup_{C_{i-1}} C_i,$$ with $p-1$ copies of $C_i$ on the left of $C_k$ and $\ell-p$ copies on the right. We can picture a minimal $X$-labeling of $B^k(i,\ell,p)$ as \[ \xymatrixcolsep{1pc}\xymatrixrowsep{1pc}\xymatrix{ & B^k(i,\ell,p) &  \\ x_1 \ar[ru] & \cdots & x_{\ell} ,\ar[lu]}\] where $x_j\in X(i)$ for $i\neq p$ and $x_p\in X(k)$ are such that $t(x_j)=s(x_{j+1})$ for $j\notin\{ p-1,p\}$, $t(x_{p-1})=s^{k-i+1}(x_p)$, and $t^{k-i+1}(x_p)=s(x_{p+1})$. Similarly, $$\widehat{B^k}(k,\ell)=C_k\cup_{C_{k-1}}\dots\cup_{C_{k-1}}C_k$$ with $\ell$ copies of $C_k$. A minimal $X$-labeling of $B^{k}(k,\ell)$ is of the form \[ \xymatrixcolsep{1pc}\xymatrixrowsep{1pc}\xymatrix{ & B^{k}(k,\ell) &  \\ x_1 \ar[ru] & \cdots & x_{\ell} ,\ar[lu]}\] with $x_j\in X(k)$ for all $j$ and $t(x_j)=s(x_{j+1})$. 
	
\end{remark}

\begin{remark}
	
	Let $x\in\Tree_n^{\OO,\EE}(X)_{k}$ be nondegenerate, with $H(x)=1$. Then $m(x)=M(x)$ and $\varphi_X(x)$ is a block, with $m(\varphi_X(x))=m(x)$. If $m(x)=k$, then $\lambda(v)=\circ_{k,k}$ for all $v\in I(T)$ and $\varphi_X(x)=B^{k}(k,\ell)$ for some $\ell$. If $m(x)=i<k$, then $\lambda(v)$ can take the values $\circ_{i,i}$, $\circ_{k,i}$ and $\circ_{i,k}$. We then have $\varphi_X(x)=B^k(i,\ell,p)$, for some $\ell$ and $p$. 
	
\end{remark}

Now we show that $\varphi_X$ is injective on $m$-ordered trees, starting with the case $H(x)=1$.

\begin{lemma}\label{mconstant}
	
	Let $x,y\in \Tree_n^{\OO,\EE}(X)_{k}$ be nondegenerate. Suppose $H(x)=H(y)=1$ and $\varphi_X(x)=\varphi_X(y)$. Then $x\stackrel{\epsilon}{=} y$.
\end{lemma}

\begin{proof}
	Let $B=\varphi_X(x)=\varphi_X(y)$. First suppose $m(x)=m(y)=i<k$, so all internal vertices in $x$ and $y$ are labeled by $\circ_{i,k}$, $\circ_{k,i},\circ_{i,i}$. Then $B=B^k(i,\ell,p)$ for some $\ell$ and $p$. Using the relations $(\circ_{k,i,i})$ and $(\circ_{i,i,k})$ one can eliminate all instances of $\circ_{i,i}$. Using the relation $(\circ_{i,k,i})$ one can move all instances of $\circ_{i,k}$ to the top of the tree and all instances of $\circ_{k,i}$ to the bottom. After applying this procedure to $x$ and $y$, we obtain trees of the form  \[\circ_{k,i}\to\cdots\to\circ_{k,i}\to\circ_{i,k}\to\cdots\to\circ_{i,k}\] (where the leaves have been omitted from the notation). The composite of these trees must be $B^k(i,\ell,p)$, so we must have $p-1$ instances of $\circ_{k,i}$ and $\ell-p$ instances of $\circ_{i,k}$. This means the trees are identical, so $x\stackrel{\epsilon}{=}y$.
	
	Now suppose $m(x)=m(y)=k$, so all internal vertices in $x$ and $y$ are labeled by $\circ_{k,k}$. The number of internal vertices is fixed by the composite $B$. So we can easily transform $x$ and $y$ into the same canonical form, by using the relation $(\circ_{k,k,k})$.\end{proof}	

Now we need some Lemmas to prepare the proof of the inductive step.

\begin{lemma}\label{H-1}
	
	Let $x\in (T_n^{\sd})^2(X)$ consist of a block $B$ labeled by nontrivial diagrams $D_i\in T_n^{\sd}(X)$ such that $M(D_i)<m(B)$ for all $i$. Then $H(D_i)=H(\mu_X^{\sd}(x))-1$ for all $i$.	
	
\end{lemma}	

\begin{proof}
	
One checks that $\{j:\ell_j(\mu_X^{\sd}(x))>1\}=\{j:\ell_j(D_i)>1\}\cup\{m(B)\}$ for all $i$.\end{proof}

\begin{lemma}\label{aux}
	
Let $x\in\Tree_n^{\OO,\EE}(\Tree_n^{\OO,\EE}(X))$ be nondegenerate, with $H(x)=1$, and suppose each of the trees $y_i$ labeling the leaves of $x$ are also nondegenerate and have $M(y_i)<m(x)$. Then $H(y_i)=H(\mu^{\OO}_X(x))-1$.	
	
\end{lemma}	

\begin{proof}
This follows easily from Lemmas \ref{mMH} and \ref{H-1}.\end{proof}

\begin{lemma}\label{cut}
	
	Let $x=(T,r,\lambda,\phi)\in\Tree_n^{\OO}(\Tree_n^{\OO,\EE}(X))_{k}$ be nondegenerate, with nondegenerate trees labeling its leaves. Suppose $H(x)=1$ and $\mu^{\OO}(x)$ is $\stackrel{\epsilon}{=}$-compatible. Then $x$ is $\stackrel{\epsilon^2}{=}$-compatible.
	
\end{lemma}		

\begin{proof}
	
	We proceed by induction on the height $h$ of $x$. If $h=0$ there is nothing to prove. If $h=1$, then the statement is also clearly true. So let $h\geq 2$. 
	
	Let $m(x)=k$ and assume both subtrees of $x$ with root at distance $1$ from $r$ have height $\geq 1$ (if this is not the case, then a similar argument applies). Then we can write \[x=\vcenter{\vbox{\xymatrixcolsep{1pc}\xymatrixrowsep{1pc}\xymatrix@1{& & \circ_{k,k} & & \\ & \circ_{k,k}\ar[ru] & & \circ_{k,k}\ar[lu] & \\ a\ar[ru] & & b\text{  }c\ar[lu]\ar[ru] & & d\ar[lu]}}}\] with $a,b,c,d\in\Tree_n^{\OO}(\Tree_n^{\OO,\EE}(X))$.  By induction, $\xymatrixcolsep{1pc}\xymatrixrowsep{1pc}\xymatrix{a\ar[r] & \circ_{k,k} & b\ar[l]}$ and $\xymatrixcolsep{1pc}\xymatrixrowsep{1pc}\xymatrix{c\ar[r] & \circ_{k,k} & d\ar[l]}$ are $\stackrel{\epsilon^2}{=}$-compatible. We need to show that $s\left(\vcenter{\vbox{\xymatrixcolsep{1pc}\xymatrixrowsep{1pc}\xymatrix@1{a\ar[r]&\circ_{k,k}&b\ar[l]}}}\right)\stackrel{\epsilon^2}{=} t\left(\vcenter{\vbox{\xymatrixcolsep{1pc}\xymatrixrowsep{1pc}\xymatrix@1{c\ar[r]&\circ_{k,k}& d\ar[l]}}}\right),$ which is equivalent to $s(b)\stackrel{\epsilon^2}{=} t(c)$. Since $\mu^{\OO}(x)$ is $\stackrel{\epsilon}{=}$-compatible, we know that $s(\mu^{\OO}(b))\stackrel{\epsilon}{=} t(\mu^{\OO}(c))$. This means that $\xymatrixcolsep{1pc}\xymatrixrowsep{1pc}\xymatrix{\mu^{\OO}(b)\ar[r] & \circ_{k,k} & \mu^{\OO}(c)\ar[l]}$ is $\stackrel{\epsilon}{=}$-compatible. By induction, this implies $\xymatrixcolsep{1pc}\xymatrixrowsep{1pc}\xymatrix{b\ar[r] & \circ_{k,k} & c\ar[l]}$ is $\stackrel{\epsilon^2}{=}$-compatible, so $s(b)\stackrel{\epsilon^2}{=} t(c)$, as desired.
	
	Now let $m(x)=i<k$. We assume that $x$ is of the form \[x=\vcenter{\vbox{\xymatrixcolsep{1pc}\xymatrixrowsep{1pc}\xymatrix@1{& & \circ_{k,i} & & \\ & \circ_{i,k}\ar[ru] & & \circ_{i,i}\ar[lu] & \\ a\ar[ru] & & b\text{  }c\ar[lu]\ar[ru] & & d\ar[lu]}}}\]  with $a,b,c,d\in\Tree_n^{\OO}(\Tree_n^{\OO,\EE}(X))$ (in other cases similar arguments apply). Now we need to show $s^{k-i+1}\left(\vcenter{\vbox{\xymatrixcolsep{1pc}\xymatrixrowsep{1pc}\xymatrix@1{a\ar[r]&\circ_{i,k}& b\ar[l]}}}\right)\stackrel{\epsilon^2}{=} t\left(\vcenter{\vbox{\xymatrixcolsep{1pc}\xymatrixrowsep{1pc}\xymatrix@1{c\ar[r]&\circ_{i,i}& d\ar[l]}}}\right),$ which is equivalent to $s^{k-i+1}(b)\stackrel{\epsilon^2}{=} t(c)$. This follows by induction as in the previous case.\end{proof}		

\begin{remark}
	
	Example \ref{grafting} shows that this Lemma would be false without the condition that $H(x)=1$.
\end{remark}

\begin{lemma}\label{decomp}
	
Every nontrivial diagram $D\in T_n^{\sd}(X)_{k}$ which is not a block has a unique decomposition of the form \[D=\vcenter{\vbox{\xymatrixcolsep{1pc}\xymatrixrowsep{1pc}\xymatrix{& B^k & \\ D_1\ar[ru] & \cdots & D_{\ell}\ar[lu]}}}\]	with $B^k$ a block, $D_i$ nontrivial, $m(B^k)=M(D)$, $M(D_i)<M(D)$ and $\ell = \ell_{M(D)}(D)$.
	
\end{lemma}		

\begin{proof}
	
If $M(D)=k$, then we must take \begin{center}$B^k=B^k(k,\ell)$ and $D_j=(v^1(D),\cdots,v^{k-1}(D),(v^k(j))).$ \end{center} If $M(D)=M<k$, then $v^{M+1}(D)=(p)$ for some $p$ and we must take \begin{center}$B^k=B^k(M,\ell,p)$ and 
$D_j = \begin{cases}
(v^1(D),\cdots,v^{M-1}(D),(v^M(j))) &\text{if $j\neq p$}\\
(v^1(D),\cdots,v^{M-1}(D),(v^M(p)),(1),\cdots,(1)) &\text{if $j=p$}.
\end{cases}$
\end{center} Since $D$ is not a block, we must have $\ell_m(D)>1$ for some $m<M$, therefore the $D_j$ are nontrivial. \end{proof}	

\begin{remark}
 
Using the Lemma above and induction on $H(D)$, one can show that every nontrivial diagram $D$ has a unique strictly ordered block decomposition. By this we mean a $B(D)\in (T_n^{\sd})^{H(D)}(X)$ whose composite is $D$ and such that every diagram appearing is a block, and moreover $m(B_1)<m(B_2)$ whenever $B_1$ labels $B_2$. We do not provide more details as we will not need this construction for the proof of our main results. 
 
\end{remark}

\begin{lemma}\label{mordered}
	
Let $x,y\in\Tree_n^{\OO,\EE}(X)$ be $m$-ordered and suppose $\varphi_X(x)=\varphi_X(y)$. Then $x\stackrel{\epsilon}{=}y$.	
	
\end{lemma}	

\begin{proof}
	
We proceed by induction on $H=H(x)=H(y)$. The base case $H=0$ is trivial, because then $x,y$ must have height zero. The case $H=1$ is Lemma \ref{mconstant}. So suppose $H>1$. Let $D=\varphi_X(x)=\varphi_X(y)$. Let $\tilde{x}\in\Tree_n^{\OO}(\Tree_n^{\OO,\EE}(X))$ with $\mu^{\OO}_X(\tilde{x})=x$ be obtained by taking the maximal subtree of $x$ with $H(\tilde{x})=1$ contanining the root. Define $\tilde{y}$ in the same way. By Lemma \ref{cut}, we have $\tilde{x},\tilde{y}\in\Tree_n^{\OO,\EE}(\Tree_n^{\OO,\EE}(X))$. Now $\varphi_X(x)$ and $\varphi_X(y)$ can be computed by applying the composite \[\xymatrix@1{\Tree_n^{\OO,\EE}(\Tree_n^{\OO,\EE}(X))\ar[rr]^-{\Tree_n^{\OO,\EE}(\varphi_X)}&& \Tree_n^{\OO,\EE}(T_n^{\sd}(X))\ar[rr]^-{\varphi_{T_n^{\sd}(X)}} && (T_n^{\sd})^2(X)\ar[r]^-{\mu_X^{\sd}} & T_n^{\sd}(X)}\] to $\tilde{x}$ and $\tilde{y}$. Since $x,y$ are $m$-ordered, the result of applying the first two maps to $\tilde{x},\tilde{y}$ is the decomposition of $D$ from Lemma \ref{decomp}, consisting of a block $B$ labeled by nontrivial diagrams $D_i$ with $m(B)=M(D)$ and $M(D_i)<m(B)$. If we let $\bar{x},\bar{y}$ be the images of $\tilde{x},\tilde{y}$ under the first map, we then have $\varphi_{T_n^{\sd}(X)}(\bar{x})=\varphi_{T_n^{\sd}(X)}(\bar{y})$ by uniqueness of the decomposition, so $\bar{x}\stackrel{\epsilon}{=}\bar{y}$ by Lemma \ref{mconstant}, since $H(\bar{x})=H(\bar{y})=1$.

 Now the leaves of $\tilde{x}$ and $\tilde{y}$ are labelled by trees $x_i,y_i$ and we have $\varphi_{X}(x_i)=D_{f(i)}$, $\varphi_{X}(y_i)=D_{g(i)}$ for some functions $f,g$ of the indices, so $\varphi_X(x_i)=\varphi_X(y_j)$ whenever $f(i)=g(j)$. We have $M(x_i)<m(\tilde{x})$ and $M(y_i)<m(\tilde{y})$ because $x,y$ are $m$-ordered, therefore $H(x_i)=H(y_i)=H-1$ by Lemma \ref{aux}. Now by induction we have $x_i\stackrel{\epsilon}{=}y_j$ whenever $f(i)=g(j)$. This allows us to lift the equivalence $\bar{x}\stackrel{\epsilon}{=}\bar{y}$ to an equivalence $\tilde{x}\stackrel{\epsilon^2}{=}\tilde{y}$, which implies $x\stackrel{\epsilon}{=}y$.\end{proof}	

\begin{proposition}\label{injective}
	
The map $\tilde{\varphi}\colon \widetilde{\Tree}_n^{\OO,\EE}\to T^{\sd}_n$ is injective.
	
\end{proposition}	

\begin{proof}
	
This follows immediately from Lemmas \ref{nondegenerate1}, \ref{nondegenerate2}, \ref{ordered} and \ref{mordered}.\end{proof}	

\subsection{$\tilde{\varphi}$ is surjective}

\begin{notation}
 
Let $D$ be a simple $k$-string diagram, $S\subset\langle\ell_{m}(D)\rangle$ and denote by $v^m_S$ the restriction of $v^m$ to $S$. If we suppose that $\im (v^{m+1})\subset S$ (or $m=k$), then

\[\xymatrixcolsep{2.5pc}\xymatrix@1{\langle\ell_k\rangle\ar[r]^-{v^k} & \cdots\ar[r]^-{v^{m+2}} & \langle\ell_{m+1}\rangle\ar[r]^-{v^{m+1}} & S \ar[r]^-{v^m_S} & \langle\ell_{m-1}\rangle \ar[r]^-{v^{m-1}} & \cdots \ar[r]^-{v^1} & \langle\ell_0\rangle=\langle 1 \rangle}\] is still a $k$-level tree with crossings, so it is uniquely isomorphic to a unique simple $k$-string diagram, which we denote by $D_S$.
 
\end{notation}

\begin{proposition}
	
The map $\tilde{\varphi}\colon \widetilde{\Tree}_n^{\OO,\EE}\to T^{\sd}_n$ is surjective.	
	
\end{proposition}

\begin{proof}
	
Let $D\in T_n^{\sd}(X)_{k}$ and let $M_0(D)=\max\{i:\ell_i(D)>0\}$. If $M_0(D)<k$ then $D=\varphi_{T_n^{\sd}(X)}(x)$, where $x=(s^{k-M_0(D)}(D)\to u_{M_0(D)-1}\to\cdots\to u_k)$. Now it is enough to show that $s^{k-M_0(D)}(D)$ is in the image of $\varphi_X$. This means we can reduce to the case where $M_0(D)=k$, i.e. $D$ is nondegenerate.

Now we do induction on $M(D)=\max\{i:\ell_i(D)>1\}$ and $\ell_{M(D)}(D)$. It is useful here to extend the definition of $M$ by setting $M(D)=0$ when $D$ is a trivial nondegenrate diagram. If $M(D)=0$, then $D=x\to I_k$, for some $x\in X_k$, and then $D=\varphi_X(x)$, where $x$ here denotes the tree with one vertex, labelled by $x$. 

Now let $D\in T_n^{\sd}(X)$ be nontrivial and $M:=M(D)$. Suppose $M<k$. Then $v^{M+1}$ is a sequence of length $1$, with values in $\langle\ell_{M}\rangle$. If $v^{M+1}(1)<\ell_{M}$, then we have $$D=\varphi_{T_n^{\sd}(X)}\left(\vcenter{\vbox{\xymatrixcolsep{1pc}\xymatrix{ (s^{k-M}(D))_{\{\ell_M\}} \ar[r] & \circ_{M,k} & D_{\langle\ell_{M}\rangle\setminus\{\ell_M\}}\ar[l]}}}\right).$$ Now $M((s^{k-M}(D))_{\{\ell_M\}})\leq M-1$ and $\ell_M(D_{\langle\ell_{M}\rangle\setminus\{\ell_M\}})=\ell_M-1,$ so by induction there exist $x,y\in\Tree_n^{\OO,\EE}(X)$ such that $\varphi_X(x)=(s^{k-M}(D))_{\{\ell_M\}}$ and $\varphi_X(y)=D_{\langle\ell_{M}\rangle\setminus\{\ell_M\}}$. Now let $z=\xymatrixcolsep{1pc}\xymatrix{ x\ar[r]& \circ_{M,k} &y\ar[l]}$. We have $s((s^{k-M}(D))_{\{\ell_M\}})=t^{k-M+1}(D_{\langle\ell_{M}\rangle\setminus\{\ell_M\}})$ and $\varphi_X$ is injective by Proposition \ref{injective}, so $s(x)\stackrel{\epsilon}{=}t^{k-M+1}(y)$, which means $z$ is $\stackrel{\epsilon}{=}$-compatible and its image under $\varphi_X$ is $D$. 

In the case where $M<k$ and $v^{M+1}(1)=\ell_M$ we apply the same argument with $D=\varphi_{T_n^{\sd}(X)}\left(\vcenter{\vbox{\xymatrixcolsep{1pc}\xymatrix{ D_{\langle\ell_{M}\rangle\setminus\{1\}} \ar[r] & \circ_{k,M} & (s^{k-M}(D))_{\{1\}\subset\langle\ell_M\rangle} \ar[l]}}}\right).$ Finally, in the case where $M=k$ we use $D=\varphi_{T_n^{\sd}(X)}\left(\vcenter{\vbox{\xymatrixcolsep{1pc}\xymatrix{D_{\{\ell_k\}} \ar[r] & \circ_{k,k} & D_{\langle\ell_{k}\rangle\setminus\{\ell_k\}}\ar[l]}}}\right).$ \end{proof}

\section{An inductive characterization of $T_n^{\sd}$-algebras}

In the previous section we gave a finite presentation of the monad $T_n^{\sd}$. This gives a characterization of $T_n^{\sd}$-algebras, which we also call $n$-sesquicategories, as $n$-globular sets equipped with binary operations of whiskering and composition and unit morphisms, satisfying associativity and unitality relations. Now we want to give another characterization of $n$-sesquicategories, which is an inductive characterization, of the same nature as the well known description of strict $n$-categories as categories enriched in $\Cat_{n-1}$.

\begin{definition}

We define a category $(\Sesq_{n-1}\Cat)^{\flat}$. Its objects are pairs $(\CC,\underline{\Hom}_\CC)$ where $\CC$ is a category and $\underline{\Hom}_\CC$ is a lift  \[\xymatrix{ & \Sesq_{n-1}\ar@{}[d]|-{=}\ar[rd]^{(-)_0} &  \\\CC^{op}\times\CC\ar@{.>}[ru]^{\underline{\Hom}_\CC}\ar[rr]_{\Hom_\CC} & & \Set .}\] A morphism $(\CC,\underline{\Hom}_\CC)\to(\DD,\underline{\Hom}_\DD)$ consists of a functor $\CC\to\DD$ and a natural transformation \[\xymatrix{\CC^{op}\times\CC\ar[rr]^{\underline{\Hom}_\CC}\ar[rd]_{F^{op}\times F} & \ar@{=>}[d] & \Sesq_{n-1} \\ & \DD^{op}\times\DD\ar[ru]_{\underline{\Hom}_\DD} & }\] such that

 \[\vcenter{\vbox{\xymatrix{\CC^{op}\times\CC\ar[rr]^{\underline{\Hom}_\CC}\ar[rd]_{F^{op}\times F} & \ar@{=>}[d] & \Sesq_{n-1} \ar@{}[d]|-{=}\ar[rd]^{(-)_0} & \\ & \DD^{op}\times\DD\ar[ru]^{\underline{\Hom}_\DD}\ar[rr]_{\Hom_{\DD}} & & \Set}}}=\vcenter{\vbox{\xymatrix{ & \Sesq_{n-1}\ar@{}[d]|-{=}\ar[rd]^{(-)_0} &  \\\CC^{op}\times\CC\ar[ru]^{\underline{\Hom}_\CC}\ar[rr]^{\Hom_\CC}\ar[rd]_{F^{op}\times F} & \ar@{=>}[d] & \Set \\ & \DD^{op}\times\DD\ar[ru]_{\Hom_\DD} & }}}.\]

\end{definition}

\begin{remark}
	
The category of $\Sesq_{n-1}$-enriched categories would be denoted by $\Sesq_{n-1}\Cat$. We use $\flat$ to indicate that the notion we are using is weaker than enrichment.	
	
\end{remark}		

\begin{theorem}\label{inductive}
	
The categories $\Sesq_1$ and $\Cat$ are isomorphic. The category $\Sesq_n$ is isomorphic to $(\Sesq_{n-1}\Cat)^{\flat}$.
	
\end{theorem}	

It's clear that $T^{\sd}_1$-algebras are the same thing as categories, so we just need to construct an isomorphism $\Sesq_n\simeq(\Sesq_{n-1}\Cat)^{\flat}$. 

\subsection{The functor $\Sesq_n\to(\Sesq_{n-1}\Cat)^{\flat}$.}

\begin{definition}
	
Given an $n$-sesquicategory $\CCC$, we denote by $\CCC_{\leq 1}$ the category whose underlying $1$-globular set is the truncation of $\CCC$ and whose $T_1$-algebra structure is the restriction of the $T_n$-algebra structure on $\CCC$. 	
	
\end{definition}	

\begin{definition}
	
Given an $n$-globular set $X$, and $x,y\in X(0)$, we denote by $\Hom_{X}(x,y)$ the $(n-1)$-globular set with $$\Hom_{X}(x,y)(k-1)=\{z\in X(k):s^{k}(z)=x,t^{k}(z)=y\}.$$ 	
	
\end{definition}	

Now we must define a $T_{n-1}^{\sd}$-algebra structure on $\Hom_{\CCC}(x,y)$ when $\CCC$ is an $n$-sesquicategory. The idea is that a $\Hom_\CCC(x,y)$-labelling of a $(k-1)$-diagram $D$ is the same thing as a $\CCC$-labeling of a certain $k$-diagram $S(D)$ obtained by shifting the dimensions of the cells up by one and introducing two new $0$-cells, which are labelled by $x$ and $y$. We now make this precise. 

\begin{definition}
	
We define a map $S\colon \sd_{n-1}(k-1)\to\sd_{n}(k)$ by $$S(v^1,\cdots,v^{k-1}):=(1,v^1,\cdots,v^{k-1}).$$ For any $D\in\sd_{n-1}(k-1)$, the $(k)$-globular set $\widehat{S(D)}$ has two $0$-cells, which we denote by $s$ and $t$, since they are equal to $s^m(z)$ and $t^{m}(z)$, respectively, for any $m\geq 1$ and $z\in\widehat{S(D)}(m)$.  	
	
\end{definition}

\begin{definition}
	
Given $D\in\sd_{n-1}(k-1)$ and $X\in\gSet_{k}$ we denote by $\gSet_{k}^{x,y}(\widehat{S(D)},X)$ the set of maps sending $s$ to $x$ and $t$ to $y$.
	
\end{definition}		

\begin{lemma}
	
Given $D\in\sd_{n-1}(k-1)$ and $X\in\gSet_{k}$, we have $$\gSet_{k-1}(\widehat{D},\Hom_X(x,y))=\gSet_k^{x,y}(\widehat{S(D)},X).$$	
	
\end{lemma}	

\begin{proof}
	
This is clear from the definitions.\end{proof}	

\begin{lemma}
	
Let $\CCC$ be an $n$-sesquicategory, $x,y\in \CCC(0)$ and $D\in\sd_{n-1}(k-1)$. The composition map $\gSet_k^{x,y}(\widehat{S(D)},\CCC)\to\CCC(k)$ factors through the inclusion $$\Hom_{\CCC}(x,y)(k-1)\hookrightarrow\CCC(k).$$
	
\end{lemma}	

\begin{proof}This is clear from the definitions.\end{proof}

\begin{definition}
	
Let $\CCC$ be an $n$-sesquicategory and $x,y\in \CCC(0)$. We define a $T_{n-1}^{\sd}$-algebra structure on $\Hom_\CCC(x,y)$ by letting $D\in\sd_{n-1}(k-1)$ act via the map $$\gSet_{k-1}(\widehat{D},\Hom_\CCC(x,y))=\gSet_k^{x,y}(\widehat{S(D)},\CCC)\to\Hom_{\CCC}(x,y)(k-1).$$	
	
\end{definition}

\begin{example}
	
Suppose $f,g,h\colon x\to y$ are $1$-morphisms in $\CCC$ and $\alpha\colon f\to g$, $\beta\colon g\to h$ are $2$-morphisms. Then the $\Hom_\CCC(x,y)$-labeled diagram on the left corresponds to the $\CCC$-labeled diagram on the right. \begin{center}\begin{tabular}{lcccr}\includegraphics[scale=1.5,align=c]{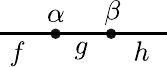} & & & & \includegraphics[scale=1.5,align=c]{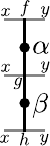} \end{tabular}\end{center} 	
	
\end{example}	

\begin{definition}
	
Let $\CCC$ be an $n$-sesquicategory. We extend $\Hom_{\CCC}$ to a functor $$\Hom_{\CCC}\colon \CCC_{\leq 1}^{op}\times\CCC_{\leq 1}\to\Sesq_{n-1}$$ as follows. For $1$-morphisms $f\colon x\to x'$ and $g\colon y\to y'$ in $\CCC_{\leq 1}$, we define $$\Hom_{\CCC}(f,g)_{k}\colon \Hom_{\CCC}(x',y)_{k}\to\Hom_{\CCC}(x,y')_{k}$$ by restricting the composite \[\xymatrixcolsep{4.5pc}\xymatrix@1{\CCC(k+1)\ar[r]^-{g\times\id\times f} & \CCC(1)\times_{\CCC(0)}\CCC(k+1)\times_{\CCC(0)}\CCC(1)\ar[r]^-{B^{k+1}(1,3,2)}&\CCC(k+1)}.\] 	
	
\end{definition}	

\begin{example}
	
When $k=1$, the map $\CCC(2)\to\CCC(2)$ from the previous definition can be depicted as$$\includegraphics[scale=1.5]{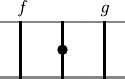} .$$
	
\end{example}		

\begin{definition}
	
Let $F\colon \CCC\to\DDD$ be a morphism in $\Sesq_n$. We define a natural transformation $$\Hom_F\colon \Hom_\CCC\Rightarrow\Hom_\DDD\circ (F_{\leq 1}^{op}\times F_{\leq 1})$$ by letting the map $\Hom_\CCC(x,y)\to\Hom_\DDD(F(x),F(y))$ be simply the restriction of the map $F\colon \CCC\to\DDD$.	
	
\end{definition}

\begin{definition}
	
	We define the functor $\Sesq_n\to(\Sesq_{n-1}\Cat)^{\flat}$ on objects by mapping $\CCC$ to the pair $(\CCC_{\leq 1},\Hom_{\CCC})$ and on morphisms by mapping $F$ to the pair $(F_{\leq 1},\Hom_F)$.	
	
\end{definition}

\subsection{The functor $(\Sesq_{n-1}\Cat)^{\flat}\to\Sesq_n$.}

\begin{definition}
	
Let $(\CC,\underline{\Hom}_\CC)$ be an object in $(\Sesq_{n-1}\Cat)^{\flat}$. We define $\CC^+\in\gSet_n$ by letting $\CC^+(0)=\obj(\CC)$ and $$\CC^+(k)=\coprod_{x,y\in\obj(\CC)}\underline{\Hom}_\CC(x,y)(k-1)$$ with obvious source and target maps. Given a morphism $(\CC,\underline{\Hom}_{\CC})\to(\DD,\underline{\Hom}_{\DD})$ in $(\Sesq_{n-1}\Cat)^{\flat}$, there is an obvious way to define a map of $n$-globular sets $\CC^+\to\DD^+$.	
	
\end{definition}

This is clearly inverse to the previously defined construction, at the level of $n$-globular sets. Now we need to define a $T_n^{\sd}$-algebra structure on $\CC^+$. We do this by saying how each generator in $\OO_n$ acts and showing that the relations $\EE_n$ are satisfied.

\begin{definition}

 We define the action of generators $\circ_{i,j}$ and $u_i$ on $\CC^+$, for $i,j>1$. A $\CC^+$-labeling of $\circ_{i,j}$, for $i,j>1$, corresponds to an $\underline{\Hom}_\CC(x,y)$-labeling of $\circ_{i-1,j-1}$, where $x,y$ are the labels of the $0$-cells in $\widehat{\circ_{i,j}}$. Similarly a $\CC^+$-labeling of $u_i$, for $i>1$, corresponds to an $\underline{\Hom}_\CC(x,y)$-labeling of $u_{i-1}$, where $x,y$ are the labels of the $0$-cells in $\widehat{u_i}$. This means we can define the action of these generators on $\CC^+$ by using the $T_{n-1}^{\sd}$-structure on $\underline{\Hom}_\CC(x,y)$.
 
\end{definition} 

The above definition guarantees that all relations involving generators $\circ_{i,j}$ and $u_i$ for $i,j>1$ are satisfied, since they follow from relations in the $T_{n-1}^{\sd}$-algebras $\underline{\Hom}_\CC(x,y)$. 

\begin{definition}
	
We now define the action of generators $\circ_{1,j}$, $\circ_{j,1}$ and $u_1$ on $\CC^+$. The action of $\circ_{1,j}$ is a map $\CC^+(1)\times_{\CC^+(0)}\CC^+(j)\to\CC^+(j)$. We have \[\begin{aligned}\CC^+(1)\times_{\CC^+(0)}\CC^+(j) & =\coprod_{x,y,y'\in\obj(\CC)}\Hom_{\CC}(y,y')\times\underline{\Hom}_{\CC}(x,y)(j-1) \\ & =\coprod_{x,f\colon y\to y'}\underline{\Hom}_{\CC}(x,y)(j-1)\end{aligned}\] and we let $\circ_{1,j}$ act on the $(x,f)$ factor by the map $$\underline{\Hom}_{\CC}(x,f)(j-1)\colon \underline{\Hom}_{\CC}(x,y)(j-1)\to\underline{\Hom}_{\CC}(x,y')(j-1)\subset\CC^+(j).$$ We define the action of $\circ_{j,1}$ in a similar way. The generator $u_1$ acts by the map $\CC^+(0)\to\CC^+(1)$ which takes an object $x\in\CC$ to the identity $1$-morphism $x\to x$. 
	
\end{definition}	

Given this definition, the remaining relations follow simply from the fact that $\underline{\Hom}_{\CC}$ is a functor. We have therefore constructed a functor $(\Sesq_{n-1}\Cat)^{\flat}\to\Sesq_n$. It is now easy to check that the pair of functors we have constructed are inverse to each other. This concludes the proof of Theorem \ref{inductive}.

\end{document}